\documentclass[a4paper, 12pt]{article}
\usepackage{enumerate, theorem}
\usepackage{amsmath, amsfonts, amssymb}
\usepackage[height=22.5cm, width=15cm]{geometry}
\usepackage[all]{xy}
\usepackage{mathrsfs, yfonts}

\DeclareMathOperator{\id}{id}

\DeclareMathOperator{\Sym}{Sym}

\DeclareMathOperator{\C}{\mathbb{C}}
\DeclareMathOperator{\codim}{codim}

\newcommand{\parag}[1]{\paragraph{\sc{#1.}} }

\newtheorem{thm}{Theorem}[subsection]
\newtheorem{defn}[thm]{Definition}
\newtheorem{cor}[thm]{Corollary}
\newtheorem{prop}[thm]{Proposition}
\newtheorem{lemma}[thm]{Lemma}

\begin{document}

\title{On the nearly smooth complex spaces}

\author{Daniel  Barlet\footnote{Institut Elie Cartan, G\'eom\`{e}trie,\newline
Universit\'e de Lorraine, CNRS UMR 7502   and  Institut Universitaire de France.} \ and J\'on  Magn\'usson\footnote{Department of Mathematics, School of Engineering and Physical Sciences, University of Iceland.}.}

\date{13/10/2017}

\maketitle

\parag{Abstract} We introduce a class of normal complex spaces having only mild singularities (close to quotient singularities) for which we generalize the notion of a (analytic) fundamental class for an analytic cycle and also the notion of a relative fundamental class for an analytic family of cycles. We also generalize to these spaces the geometric intersection theory for analytic cycles  with rational positive coefficients and show that it behaves well with respect to analytic families of cycles.\\
We prove that this  intersection theory has most of the usual properties of the standard geometric intersection theory on complex manifolds, but with the exception that the intersection cycle of two cycles with positive integral coefficients that intersect properly may have rational coefficients.

\parag{AMS classification} 32 C 20- 32 C 25- 32 C 36.

\parag{Key words} Quotient Singularity- The Sheaf $\omega_{X}^{\bullet}$ - Fundamental Class of Cycles- Analytic Family of Cycles -Geometric Intersection Theory.

\tableofcontents

\section*{Introduction}

The theory of intersection for complex analytic cycles (with positive integral coefficients) on a complex manifold is more or less standard; see for instance [D. 69] or [T. 95] (see also  [Fult] or [FMP] for the intersection theory in algebraic geometry). This theory  for a complex manifold is thoroughly introduced in chapter VII of [B-M 2] where it is described in terms of fundamental classes for analytic cycles (see [B.80]) and their cup products. The aim of the present paper is to extend this theory to a certain class of normal complex spaces having ``nice enough'' singularities that  we have chosen to call  {\em nearly smooth}. Every  complex space with quotient singularities\footnote{This means that the space is locally isomorphic to a quotient of a complex manifold by a finite group of automorphisms.} is nearly smooth and to our knowledge there are no examples of a nearly smooth complex space that does not have quotient singularities. In fact this is related to a ``classical'' conjecture (see [K. 09] p. 83). But we will show that the nearly smooth complex spaces have some interesting stability properties  and it is not clear that complex spaces with quotient singularities have all of them (see for instance theorem \ref{stab.}). 

The fundamental class of an analytic cycle \ $X$ \ of codimension \ $p$ \ in a complex manifold \ $M$ \ is the global section of the local cohomology sheaf \ $\underline{H}_{|X|}^p(\Omega_M^p)$ \ defined by the integration current of \ $X$. Moreover, if two analytic cycles intersect properly on \ $M$ \ then the cup product of their fundamental classes is the fundamental class of their intersection cycle. In order to generalize this theory to the nearly smooth case we have replaced \ $\Omega_M^p$ \ by \ $\omega_M^p$, the sheaf of \ $\bar{\partial}-$closed $(p, 0)$  currents (modulo torsion) introduced in [B.78]. But since \ $M$ \ is normal  \ $\omega_M^p$ \ is the natural extension to \ $M$ \ of the sheaf of holomorphic \ $p-$forms on the non-singular part of \ $M$.

We prove that this generalized intersection theory has most of the usual properties of the standard intersection theory on manifolds, as compatibility with pull-back, projection formula etc...,  with the important exception that two cycles which have positive integral coefficients and intersect properly may have (positive) rational but not integral coefficients.

\section{Nearly-smooth complex spaces}

In this article all complexe spaces are {\em reduced} by hypothesis and analytic subsets are supposed to be {\em closed} unless otherwise stated.

\subsection{Definition and first properties.}

\medskip
Let us recall two notions that will be of importance in the paper.

\smallskip
The sheaf of \ $\bar{\partial}-$closed currents\footnote{For the notion of a {\em current} on a reduced complex space see for example [B-M 1]} of bidegree \ $(p,0)$ \ (modulo torsion) on a reduced  complex space \ $M$ \ will be denoted \ $\omega_{M}^p$. If \ $j : M\setminus S \to M$ \ is the inclusion of the smooth part of \ $M$ \ into \ $M$, then we have a natural (injective) \ $\mathcal{O}_M-$linear morphism 
$$ 
\omega_{M}^{\bullet}  \longrightarrow   j_{*}j^{*}(\Omega_{M}^{\bullet})
$$
and in the case when \ $M$  \ is normal this morphism is an isomorphism. In fact this is true, more generally, when the codimension of $S$ in $M$ is at least $2$ (see [B.78]). In this case we have a natural exterior product on the graded algebra $\bigoplus_{p\geq 0}\omega^{p}$.

\smallskip
We will denote by 
$$
h_M^{\bullet} : \Omega_{M}^{\bullet} \longrightarrow \omega_{M}^{\bullet}
$$ 
the \ $\mathcal{O}_M-$linear morphism which associates  to a homorphic $p-$form the corresponding $\bar \partial-$closed $(p, 0)-$current on $M$.

\medskip
In the sequel we shall use the following notion.

\begin{defn}\label{g-flat}
Let \ $f : M\rightarrow N$ \ be a holomorphic map between pure dimensional complex spaces. We say that $f$ is {\bf\em geometrically flat}  if  there exists an analytic family $(F_{y})_{y \in N}$ of cycles \ in $M$ \ parametrized by $N$ such that for every $y \in N$ the support $|F_{y}|$ of the cycle $F_{y}$ is the set $f^{-1}(y)$ and for $y$ very general in $N$ the cycle $F_{y}$ is reduced.\footnote{See theorem IV.3.3.1 in [B-M.1]}
\end{defn}

Roughly this means that, putting suitable non trivial  multiplicities on  the irreducible components of some fibres of \ $f$ \ over a relatively small subset in \ $N$,  the family of  fibres will form an analytic family of cycles parametrized by \ $N$

\parag{Remarks}
Let $f : M\rightarrow N$ be a holomorphic map between two pure dimensional complex spaces and put $d := \dim M - \dim N$.
\begin{enumerate}[i)]
\item
If the map \ $f$ \ is geometrically flat then it is open and {\em equidimensional}, i.e. $\dim_xf^{-1}(f(x)) = d$ for all\ $x$.
We also have a holomorphic {\em fibre map}\footnote{{\em Holomorphic} simply means that the family of cycles $(\varphi(y))_{y \in N}$ is analytic; see [B-M 1] ch.IV for a definition.}
$$
\varphi : N \to \mathcal{C}_{d}^{loc}(M)
$$
classifying the fibres of $f$ (see [B-M 1]). This is specially interesting when $f$ is proper because in this case  $\varphi$ induces a holomorphic map into the space $\mathcal{C}_{d}(M)$ of compact $d-$dimensional cycles in $M$ which is a reduced complex space (see [B.75] or [B-M 1] and [B-M 2]).\\
For instance, when $d = 0$ and $N$ is connected we obtain a holomorphic fibre map with values in $\Sym^{k}(M) := M^{k}\big/\frak{S}_{k}$ where $k$ is the degree of $q$. This is the most important case for this paper.
\item 
In the case where $N$ is a normal complex space the map $f$ is geometrically flat   if and only if it is equidimensional (see [B-M 1] ch.IV th. 3.4.1).
\end{enumerate}

\begin{lemma}\label{N-Smooth -1}
Let $q : \tilde{M} \to M$ be a holomorphic map which is proper,  finite and geometrically flat from a  normal complex space  $\tilde{M}$ to a connected reduced complex space $M$. Then $M$ is normal.
\end{lemma}

\parag{Proof} Let $f$ be a locally bounded meromorphic function on an open set $M'$ of $M$. Then $g := q^{*}(f)$ has to be holomorphic on $q^{-1}(M')$ as $\tilde{M}$ is normal and $q$ is proper.\\
Now let $\varphi : M \to  \Sym^{k}(\tilde{M})$  be the holomorphic fibre map of $q$. Put $\tilde{M}' := q^{-1}(M')$ and compose the restriction of \ $\varphi$ \ to $M'$ with the holomorphic map $\Sym^{k}(\tilde{M}') \to \C$,  which is induced by the first symetric function composed with the holomorphic map
$$
\Sym^{k}(g) : \Sym^{k}(\tilde{M}') \to \Sym^{k}(\C).
$$
Then we obtain a holomorphic map on $M'$  which is easily seen to be equal to $k.f$ at least at the generic points in $M'$. So $f$ is in fact holomorphic on $M'$ proving that  $M $ is normal.
$\hfill \blacksquare$\\

Let  $q : \tilde{M} \rightarrow M$ be a surjective holomorphic map between irreducible complex spaces and assume  that $q$  is proper and  generically  finite . Then by taking direct images of currents by \ $q$ \ we obtain a graded \ $\mathcal{O}_M-$linear map
$$
Trace_{q} : q_{*}(\omega_{\tilde{M}}^{\bullet}) \longrightarrow \omega_{M}^{\bullet}
$$
and a commutative diagram
$$
\xymatrix{
q_{*}(\omega_{\tilde{M}}^{\bullet})  \ar[rr]^{\ \ Trace_{q}} && \ \omega_{M}^{\bullet}\\ 
q_{*}(\Omega_{\tilde{M}}^{\bullet}) \ar[u]^{q_*(h_{\tilde M}^{\bullet})} && \ \Omega_{M}^{\bullet}\ar[u]_{k.h_M^{\bullet}}\ar[ll]_{q^*}} 
$$
where \ $k$ denotes the generic degree of \ $q$. Such maps will be called {\em trace maps}. In what follows we will occasionally also  use the symbol $Trace_{q}$ for the composed map $q_{*}(\Omega_{\tilde{M}}^{\bullet}) \rightarrow\omega_{M}^{\bullet}$. 

\smallskip
Now suppose moreover that \ $\tilde M$ \ is smooth, $M$ \ is normal and \ $q$ \ is finite, then we get a graded $\mathcal{O}_M-$linear pull-back morphism
$$
\hat{q}^{*} :  \omega_{M}^{\bullet} \to q_{*}(\Omega_{\tilde{M}}^{\bullet})
$$
which is compatible with the usual one, i.e. $\hat{q}^{*}\circ h_M^{\bullet} = q^*$. It is defined in the following way:\\
A section  $\sigma$ on an open set $M'$ in $M$ of the sheaf $\omega_{M}^{p}$ defines a $p-$holomorphic form on $M' \setminus S$ where $S$  is the singular set of $M$. So its pull-back by $q$ defines a holomorphic $p-$form on $q^{-1}(M')\setminus q^{-1}(S)$. But, as $M$ is normal, $S$ has codimension at least  $2$ and then, $q^{-1}(S)$ has  also codimension at least $2$. It follows that  the holomorphic $p-$form $q^{*}(\sigma)$ extends to a section of the sheaf $q_{*}(\Omega_{\tilde{M}}^{p})$ on $M'$ as $\tilde{M}$ is smooth. It is clear that this morphism is graded and  $\mathcal{O}_{M}-$linear.\\ 
   
The following proposition is the initial point that leads us to introduce the notion of a nearly smooth complex space.

\begin{prop}\label{N-Smooth 0}
Let $q : \tilde{M} \to M$ be a surjective holomorphic map which is proper,  finite and geometrically flat from a connected complex manifold $\tilde{M}$ to a complex space $M$. Then the (graded) $\mathcal{O}_{M}-$linear map $\hat{q}^{*} :  \omega_{M}^{\bullet} \rightarrow q_{*}(\Omega_{\tilde{M}}^{\bullet})$ induces an isomorphism onto a direct factor of $q_{*}(\Omega_{\tilde{M}}^{\bullet})$.
\end{prop}

\parag{Proof}
Let \ $k$ \ denote the degree of \ $q$. Then, outside of the branch locus of \ $q$, we clearly have\\
\begin{equation*}
Trace_{q}\circ \hat{q}^{*} = k.\id_{\omega_{M}^{\bullet}}  
\end{equation*}
and consequently everywhere on \ $M$, as the sheaf $\omega_{M}^{\bullet}$ has no torsion.\hfill $\blacksquare$

\parag{Remarks} 
\begin{enumerate}[i)]
\item
The proof of lemma \ref{N-Smooth -1} consists in fact of showing that the $\mathcal{O}_{M}-$linear map
$Trace_q : q_{*}(\mathcal{O}_{\tilde{M}}) \rightarrow \omega_{M}^0$ \ takes its values in \ $\mathcal{O}_M$ and hence $\mathcal{O}_M = \omega_{M}^0$.
\item
The morphism $\hat{q}$ preserves the exterior product.
\end{enumerate}

\begin{defn}\label{N-Smooth 1}
A complex space $M$  is called {\bf\em nearly smooth} if it is normal  and if  for every point \ $x$  in $M$ there exists an open neighbourhood $M'$ of \ $x$ \ and a proper, finite and surjective holomorphic map $q :\tilde{M}' \to M'$ from a connected complex manifold $\tilde{M}'$ .
\end{defn}

For instance, any complex space having only quotient singularities is a nearly smooth complex space. 

\parag{Terminology} When we have an open set  $M'$ in a nearly smooth complex space $M$ and a proper, finite and surjective holomorphic map $q : \tilde{M'} \to M'$ from a  complex manifold $\tilde{M'}$, we shall say that  $q : \tilde{M'} \to M'$  is a {\bf local model} for the nearly smooth complex space $M$.
$\hfill \square$\\

Our next proposition gives some basic properties of nearly smooth complex spaces.

\begin{prop}\label{N-Smooth 2}
A nearly smooth complex space $M$ has the following properties:
\begin{enumerate}[i)]
\item 
The de Rham complex $(\omega_{M}^{\bullet}, d^{\bullet})$ is a resolution of the constant sheaf \ $\underline{\C}$ on $M$.
\item 
$M$ is  Cohen-Macaulay.
\item For any  analytic subsets  $X$ and $Y$ in $M$ we have
$$
\codim_{M}(X\cap Y) \leq \codim_{M} X + \codim_{M} Y
$$
\item  
Any Weil divisor in $M$ is locally $\mathbb{Q}-$Cartier. So $M$ is locally  $\mathbb{Q}-$Gorenstein.
\item 
Let \ $X$ \ be an analytic subset of pure  codimension $p$  in $M$. Then $\underline{H}_{X}^{j}(\omega_{M}^{\bullet}) = 0$ for every $j \not= p$  and the sheaf $\underline{H}_{X}^{p}(\omega_{M}^{\bullet})$ has no non zero section with support in a nowhere dense  analytic subset  in $X$.
\item 
For any proper and finite holomorphic map $f : M \to N$ of a  nearly smooth complex space $M$ to a complex space $N$ the normalization of $f(M)$ is a nearly smooth complex space.
\end{enumerate}
\end{prop}

The property vi) will be widely generalized (see the theorem \ref{stab.} below).

\parag{Proof} All these properties are local so we may assume that we have a holomorphic surjective map $q : \tilde{M} \to M$ which is proper and finite of degree $k$ where $\tilde{M}$ is a connected complex manifold. \\
-- In order to prove i)  let $\alpha$ be a section on an open neighbourhood $U$ of $x \in M$ of the sheaf $\omega_{M}^{p}$ with $p \geq 1$, such that $d\alpha = 0$. Take an open neighbourhood $V \subset q^{-1}(U)$ of $q^{-1}(x)$ with Stein contractible connected components. Then there exists a holomorphic $(p-1)-$form on $V$ satisfying $d\beta = q^{*}(\alpha)$ on $V$ by the holomorphic de Rham lemma. Then $\gamma := Trace_{q}(\beta)$ is a section of the sheaf $\omega_{M}^{p-1}$  on an open neighbourhood $U' $ of $x$ in $M$ and $d\gamma = Trace_{q}(q^{*}(\alpha)) = k.\alpha$. This gives the exactness of the complex $(\omega_{M}^{\bullet}, d^{\bullet})$ in positive degrees. For $p = 0$, $q^{*}(\alpha)$ is a locally constant function on $V$ and so $Trace_{q}(q^{*}(\alpha))$ is also locally constant near $x$. So i) is proved.

\smallskip
-- As the property  ii) is local we may assume that  we dispose of a local parametrization $\pi' : M \to U$ of $M$, where $U$ is an open polydisc in $\C^{m}, m := \dim M$, where the map $\pi'$ is holomorphic proper and finite of degree $l$. Put $\pi := \pi' \circ q$. Since both \ $\tilde M$ and \ $U$ \ are smooth the direct image sheaf $\pi_{*}(\mathcal{O}_{\tilde{M}})$ is a free $\mathcal{O}_{U}-$module of rank $k.l\, ( = deg\, \pi)$. Then according to proposition \ref{N-Smooth 0} we have $q_{*}(\mathcal{O}_{\tilde{M}}) \simeq \mathcal{O}_{M}\oplus \mathcal{K}$ where $\mathcal{K}$ is the kernel of the trace map $Trace_{q}: q_{*}(\mathcal{O}_{\tilde{M}}) \to \omega_{M}^{0} = \mathcal{O}_{M}$ and consequently 
$$
(\pi')_{*}(q_{*}(\mathcal{O}_{\tilde{M}})) \simeq (\pi')_{*}(\mathcal{O}_{M}) \oplus (\pi')_{*}(\mathcal{K}) \simeq \mathcal{O}_{U}^{k.l}. 
$$
This implies that the $\mathcal{O}_{U}-$module  $(\pi')_{*}(\mathcal{O}_{M})$ is isomorphic to a direct factor of $\mathcal{O}_{U}^{k.l}$; so it is a free $\mathcal{O}_{U}-$module of finite rank (equal to $l$). Hence $M$ is Cohen-Macaulay.

\smallskip
-- Proof of iii). The analytic subsets $q^{-1}(X)$ and $q^{-1}(Y)$ in $\tilde{M}$ satisfy 
$$ 
\codim\,  \big(q^{-1}(X) \cap q^{-1}(Y)\big) \leq  \codim\,  q^{-1}(X) +  \codim\,  q^{-1}(Y)
$$
and for any analytic subset $Z$ in $M$ we have $\codim_{M}\, Z = \codim_{\tilde{M}} \, q^{-1}(Z)$. Then the equality $q^{-1}(X) \cap q^{-1}(Y) = q^{-1}(X \cap Y)$ gives  iii).

\smallskip
-- Proof of iv). It is enough to prove the statement for an effective Weil divisor $D$.  Then $q^{-1}(D)$ is a locally Cartier divisor on $\tilde{M}$ and, since the problem is local, we may assume that there exists a holomorphic function $g : \tilde{M} \to \C$ such that $q^{-1}(D) = g^{-1}(0)$. \\
Let $\varphi : M \to \Sym^{k}(\tilde{M})$ be the holomorphic fibre map of \ $q$ and let $\gamma : M \to \C$ be the composition of $\Sym^{k}(g)\circ \varphi$ with the last elementary  symetric function (the product)  $s_{k} : \Sym^{k}(\C) \to \C$. Then the Cartier divisor $\gamma^{-1}(0)$ in $M$ is equal to $k.D$, proving iv).

\smallskip
-- Proof of v).  As the analytic subset $q^{-1}(X)$ has pure codimension $p$ in the complex manifold $\tilde{M}$ the sheaf $\underline{H}_{q^{-1}(X)}^{j}(\Omega_{\tilde{M}}^{\bullet})$ vanishes for $j \not= p$ and the sheaf $\underline{H}_{q^{-1}(X)}^{p}(\Omega_{\tilde{M}}^{\bullet})$ has no non zero section with support in a nowhere dense analytic subset $Y$ \ of \ $q^{-1}(X)$. Now the functor $q_{*}$ is exact so we get 
$$
q_{*}\big(\underline{H}_{q^{-1}(X)}^{j}(\Omega_{\tilde{M}}^{\bullet})\big) \simeq \underline{H}_{X}^{j}\big(q_{*}(\Omega_{\tilde{M}}^{\bullet})\big),\quad\forall j \geq 0.
$$
and the fact that $\omega_{M}^{\bullet}$ is isomorphic to a direct factor of the $\mathcal{O}_{M}-$module $q_{*}(\Omega_{\tilde{M}}^{\bullet})$ allows us to  conclude.

\smallskip
-- Finally let us prove vi). As $f$ is proper, Remmert's direct image theorem says that $f(M)$ is an analytic subset in $N$. The induced map $g : M \rightarrow f(M)$ is surjective proper and finite, but possibly not geometrically flat. Let $\nu : Q \to f(M)$ be the normalization  of $f(M)$. Since $M$ is normal the map \ $g$ \ admits  a holomorphic lifting $\tilde{g} : M \rightarrow Q$  that is surjective, proper, finite and geometrically flat. So, if $q : \tilde{M}' \rightarrow M'$ is a local model for $M$, then the  map $\tilde{M} \rightarrow \tilde{g}(M')$ induced by \ $\tilde{g}\circ q$ \ is a local model for $Q$. Hence $Q$ is nearly smooth.
$\hfill \blacksquare$\\

The following interesting stability property of the nearly smooth complex spaces  is a wide generalization of the point vi) of the proposition \ref{N-Smooth 2}.

\begin{thm}\label{stab.}
Let $f : M \to N$ be a surjective geometrically flat holomorphic map between a nearly smooth complex space $M$  and an irreducible  complex space $N$. Then $N$ is nearly smooth.
\end{thm}

Note that if  we assume ``a priori'' that $N$ is normal, it is enough to assume that $f$ is surjective and equidimensional to conclude that $N$ is nearly smooth.\\

Remark also that no properness assumption is made on $f$.

\parag{Proof} The assertion is local on $N$, and also on $M$ because a geometrically flat map is open. So we may assume that we are in a model situation $q : \tilde{M} \to M$. Then, as the map $f\circ q$ is again geometrically flat, it is enough to prove the theorem for $M$ smooth.\\
Fix a point $y_{0}$ in $N$ and choose a generic (smooth) point $x_{0}$ in the set theoretic  fiber $f^{-1}(y_{0})$. Choose a smooth submanifold $\tilde{N}$ in an open neighbourhood of $x_{0}$ that is transversal to $f^{-1}(y_{0})$ at  $x_{0}$. Let $n$ and $d$ denote the dimensions of $\tilde{N}$ and $f^{-1}(y_{0})$. Hence there exists a $d-$scale $E := (U, B, j)$ adapted to $f^{-1}(y_{0})$,  where $j $ is an isomorphism of an open  neighbourhood of $x_{0}$ into a open neighbourhood of $\bar U \times \bar B$ in $\C^{d}\times \C^{n}$ and such that $j(x_{0}) = (0, 0)$, $j(f^{-1}(y_{0}))\cap (U\times B) = U \times \{0\}$ and  $j(\tilde{N}) \cap (U \times B) = \{0\}\times B$. Then, the geometric flatness of $f$ implies that there exists an open neighbourhood $N_{0}$ of $y_{0}$ and a  holomorphic map $ \varphi : N_{0}\times U \to \Sym^{k}(B)$ classifying the fibers of $f$ in the scale $E$. Consequently, shrinking $\tilde{N}$ if necessary, the restriction of $f$ to $\tilde{N}$ induces a geometrically holomorphic map  $\tilde{N}\rightarrow N_{0}$ which is proper and finite of degree $k$, in other words a local model for $N$.
$\hfill \blacksquare$

\subsection{Pull-back morphisms on nearly smooth spaces}

Let $q : T\rightarrow S$ be a surjective, proper, finite and geometrically flat map between complex spaces. We say that an $\mathcal{O}_S-$morphism $\hat{q}^* : \omega_S^{\bullet}\rightarrow q_*\Omega_T^{\bullet}$ is a {\bf pull-back morphism} (for $q$) if the diagram 
$$
\xymatrix{
\Omega_S^{\bullet} \ar[d]_{h_S^{\bullet}} \ar[r]^{q^*} & q_*\Omega_T^{\bullet}\\ 
\omega_S^{\bullet} \ar[ru]_{\hat{q}^*} & {}
} 
$$
is commutative. We noticed in paragraph 1.1 that the map $q$ admits such a pull-back morphism if $T$ is a complex manifold, but this is also true if the map $q$ factorizes through a complex manifold. In fact if $q = q_1\circ q_2$ where both $q_1$ and $q_2$ are surjective, proper, finite and geometrically flat maps and $q_1$ admits such a pull-back morphism then so does $q$.

\smallskip
Consider a cartesian diagram of complex spaces
$$
\xymatrix{
\tilde{Z} \ar[d]_{\pi} \ar[r]^{\ \ \tilde{f}} & T \ar[d]^{q}\\ 
Z \ar[r]^{f} & S
} 
$$
where $q$ is surjective, proper, finite and geometrically flat and $Z$ is irreducible.\footnote{The space $\tilde{Z}$ is pure dimensional, and the map $\pi$ is proper, finite and surjective but not geometrically flat in general.} Suppose also that we have an analytic cycle $G$ whose support is $\tilde{Z}$, i.e. 
$$G  = c_1\tilde{Z}_1 + \cdots + \tilde{Z}_l$$
 where $\tilde{Z}_1,\ldots,\tilde{Z}_l$ are the irreducible components of $\tilde{Z}$ and $c_1,\ldots,c_l\in\mathbb{N}^*$.
Let $[G]$ denote the {\em integration current} of $G$.  For every open subset $V$ of $Z$ and every $\gamma$ in $\Omega_{\tilde{Z}}^p(\pi^{-1}(V))$ the support of the current $\gamma\wedge[G]$ (on $\pi^{-1}(V)$) is proper over $V$ and its dircet image $\pi_*(\gamma\wedge[G])$ is an element of $\omega_Z^p(V)$. Hence we get a graded $\mathcal{O}_Z-$morphism 
$$
Trace_{G/Z} : \pi_*\Omega_{\tilde Z}^{\bullet}\longrightarrow \omega_Z^{\bullet}.
$$
If we let $d_i$ denote the degree of $q_{|\tilde{Z}_i} : \tilde{Z}_i\rightarrow Z$ and put $d := c_1d_1+\cdots +c_ld_l$ then for every $\beta$ in $\Omega_Z^{\bullet}(V)$ we get 
$$
\frac 1dTrace_{G/Z}(\pi^*\beta) = h_Z^{\bullet}(\beta).
$$
The integer $d$ will be called the degree of $G$ and denoted $\deg G$.

\smallskip
For every open subset $U$ in $S$ and every $\beta$ in $\Omega_T^{\bullet}(q^{-1}(U))$ we have 
$$
\frac 1dTrace_{G/Z}(\tilde{f}^*\beta) \ \in \ \omega_Z^{\bullet}(\pi(\tilde{f}^{-1}(q^{-1}(U))) \ = \ \omega_Z^{\bullet}(f^{-1}(U))
$$
and hence we get a graded $\mathcal{O}_Z-$morphism $\tilde{f}^*_G : q_*\Omega_T^{\bullet}\rightarrow f_*\omega_Z^{\bullet}$.

\smallskip
The following result is obvious.

\begin{lemma}\label{pull-back.-1}
In the situation described above suppose that $q$ admits a pull-back morphism $\hat{q}^* : \omega_S^{\bullet}\rightarrow q_*\Omega_T^{\bullet}$. Then the $\mathcal{O}_S-$morphism $\hat{f}_G^* := \tilde{f}^*_G\circ\hat{q}^*$ makes the diagram
$$
\xymatrix{
\Omega_S^{\bullet} \ar[d]_{h_S^{\bullet}} \ar[r]^{f^*} & f_*\Omega_Z^{\bullet}\ar[d]^{f_*(h_Z^{\bullet})}\\ 
\omega_S^{\bullet} \ar[r]^{\hat{f}_G^*} & f_*\omega_Z^{\bullet} 
} 
$$
commutative.
\hfill{$\blacksquare$}
\end{lemma}

The cycle structure on $Z\times_ST$ can come up in several different ways but two are the most important in the sequel. \\ 
--The first one is the reduced cycle given by the analytic set $Z\times_ST$.\\ 
--The other is obtained by the following procedure: assume that $S$ is connected and let $k := \deg q$. By composing  the fibre map $\varphi: S \rightarrow \Sym^{k}(T)$ of $q$ with $f$ we obtain a holomorphic map $\varphi\circ f : Z\rightarrow\Sym^k(T)$ \ that defines an analytic family \ $(X_z)_{z\in Z}$ \ of  \ $0-$cycles in \ $T$ \ parametrized by \ $Z$. Its graph cycle is an analytic cycle in \ $Z\times T$ and it will be denoted $f^*G$ if $G$ denotes the graph cycle of $q$.  Obviously \ $|f^*G| = Z\times_ST$, but  \ $f^*G$ \ is not a reduced cycle if \ $f(Z)$ \ is contained in the branched locus of \ $q$. \\
Remark that the direct image of the cycle  $f^{*}G$ in $Z$  is always  equal to $k.Z$ in this second case.

\begin{lemma}\label{pull-back.3}
Let $S$ be an irreducible complex space and $q : T \rightarrow S$ be a finite, proper and geometrically flat map. Let $Z$ be a locally closed irreducible analytic subset of $S$,  $j : Z \rightarrow S$ be the natural injection and $G$ be the graph cycle of $q$. Consider the cartesian square
$$
\xymatrix{
q^{-1}(Z) \ar[r]^{\quad \tilde{j}} \ar[d]_{\tilde{q}} & T \ar[d]^{q}\\ 
Z \ar[r]^{j} & {S}} 
$$
and suppose we have a global section $\beta$ of \ $\Omega_{T}^{p}$ such that $Trace_{q}(\beta) = 0$. Then $Trace_{j^*G/Z}(\tilde{j}^*\beta) = 0$ in $\omega_{Z}^{p}(Z)$.
\end{lemma}

\parag{Proof} Since $\omega_{Z}^{p}$ is a torsion free $\mathcal{O}_{Z}-$module the result is clear if $Z$ is not contained in the singular part of  $S$. Hence we may assume that $Z$ is smooth and contained in the singular part of  $S$. The problem being local on $S$ along $Z$ we may also assume that we have a finite and proper holomorphic map (local parametrization) $\pi : S \to U = V\times W$ where $V$ and $W$ are polydiscs centered at the origins in $\C^{r}$ and $\C^{n-r}$ with $r := \dim Z$ and $n := \dim S$, having the following properties:
\begin{itemize}
\item 
$\pi(Z) = V \times \{0\}$.
\item 
$\pi^{*}(V \times \{0\}) = k.Z$ as cycles in $S$, where $k := \deg \pi$.
\end{itemize} 

For each $w \in W$ put $Z_{w} :=  \pi^{*}(V\times\{w\})$ and $\tilde{Z}_{w} :=  q^{*}(Z_{w})$. Hence we obtain two  analytic families of $r-$cycles in $S$ and $T$ parametrized by $W$.\\
Take any $\varphi$ in $\mathscr{C}_{c}^{\infty}(Z)^{(n-r,n)} \simeq \mathscr{C}_{c}^{\infty}(V\times \{0\})^{(n-r,n)}$ and choose $\Phi_{0}$ in $\mathscr{C}^{\infty}(U)^{(n-r,n)} $ with a $W-$proper support that induces $\varphi$ on $V \times \{0\}$. Let $\Phi$ and $\tilde{\Phi}$ denote the pull-backs of $\Phi_{0}$ on $S$ and $T$ respectively. Then we get
$$
k.\left\langle Trace_{j^*G/Z}(\tilde{j}^*\beta), \varphi \right\rangle = \int_{\tilde{Z}_{0}}  \beta\wedge \tilde{\Phi}.
$$
Since  $(\tilde{Z}_{w})_{w \in W}$ is an  analytic family  and $\tilde{\Phi}$ has a $W-$propre support the function
$$
\theta(w) := \int_{\tilde{Z}_{w}}\,\beta\wedge \tilde{\Phi} = \int_{Z_w}\,Trace_q(\beta)\wedge \Phi
$$
is continuous on $W$. For all $w$ in an open dense subset $W'$ of $W$ the cycle $Z_w$ is not contained in the singular part of $S$ and consequently $\theta(w) = 0$ for all $w\in W'$. Hence the function $\theta$ is identically zero, in particular $\theta(0) = 0$, and we can conclude that $Trace_{j^*G/Z}(\tilde{j}^*\beta) = 0$.
$\hfill \blacksquare$

\begin{lemma}\label{pull-back.4}
Let $S$ be an irreducible complex space, $q : T \rightarrow S$ a finite proper and geometrically flat map and $G$ be the graph cycle of $q$. Let
$f : Z\rightarrow S$ be a holomorphic map from an irreducible complex space and denote by $\tilde{f} : Z\times_ST \rightarrow T$ and  $\pi : Z\times_ST \rightarrow Z$ the natural maps. Let $\beta$ be a global section of $\Omega_{T}^{p}$ such that $Trace_{q}(\beta) = 0$.  Then $Trace_{f^*G/Z}(\tilde{f}^*\beta) = 0$ in $\omega_{Z}^{p}(Z)$.
\end{lemma}
\parag{Proof}
Consider a stratification \ $\cdots \subseteq S_2\subseteq S_1 \subseteq S_0 = S$ \ having the following properties:
\begin{itemize} 
\item
$S_{j+1}$ \ is a nowhere dense analytic subset of $S_j$  such that $S_j\setminus S_{j+1}$ is smooth.
\item
The induced map $q^{-1}(S_j\setminus S_{j+1})\rightarrow S_j\setminus S_{j+1}$ is a covering map.
\end{itemize} 
Then there is a unique $j_0$ such that $f(Z)$ is contained in $S_{j_0}$ but not contained in $S_{j_0+1}$. Pick an irreducible component $\Sigma$ of $S_{j_0}$ that contains $f(Z)$ and let $\iota : \Sigma\rightarrow S$ denote the natural injection. Let \ $g : Z\rightarrow \Sigma$ denote the induced map and consider the commutative diagram
$$
\xymatrix{
Z\times_ST  \ar@/^2pc/[rr]^{\tilde{f}}\ar[r]^{\tilde{g}}\ar[d]_{\pi}& q^{-1}(\Sigma)\ar[d]^{q_1}\ar[r]^{\tilde{\iota}} & T\ar[d]^{q}\\ 
Z  \ar@/_2pc/[rr]^f\ar[r]^{g} & \Sigma\ar[r]^{\iota} & S
}
$$
where $\tilde{g}$, $q_1$ and $\tilde{\iota}$ are the natural maps.
Since $Z$ is irreducible and \ $\omega_Z$ \ is torsion free it is enough to show that the section $Trace_{f^*G/Z}(\tilde{f}^*\beta)$ vanishes on some open non empty subset of $Z$. Let $V$ be an open non empty $q_1-$trivializing subset of $\Sigma\setminus S(\Sigma)$ and let $W$ be an open non empty $\pi-$trivializing subset of $g^{-1}(V)$. If we let $k$ denote the degree of $q$ then there exist (not necessarily different) sections $\sigma_1,\ldots,\sigma_k$ of $q_1$ such that 
$$
Trace_{\iota^*G/\Sigma}(\tilde\iota^*\beta) = \sum_{\nu = 1}^k\sigma_{\nu}^*(\tilde{\iota}^*\beta).
$$
From lemma \ref{pull-back.3} we know that $Trace_{\iota^*G/\Sigma}(\tilde\iota^*\beta) = 0$ so we finally get
$$
Trace_{f^*G/Z}(\tilde{f}^*\beta) = \sum_{\nu = 1}^k(\sigma_{\nu}\circ g)^*(\tilde{\iota}^*(\beta)) 
= g^*\left(\sum_{\nu = 1}^k\sigma_{\nu}^*(\tilde{\iota}^*(\beta))\right) = 0.
$$
\hfill{$\blacksquare$}

\parag{Remark} 
In the situation of lemma \ref{pull-back.4} suppose that $q$ admits a pull-back morphism $\hat{q}^* : \omega_S\rightarrow q_*\Omega_T$ and let $k$ denote the degree of $q$. Then for any  $\alpha\in\omega_S^{p}(S)$ we have 
$$
Trace_{f^*G/Z}(\tilde{f}^*(\beta)) =\frac 1kTrace_{f^*G/Z}(\tilde{f}^*(\hat{q}^*(\alpha)))
$$
for any $\beta\in\Omega_{T}^{p}(T)$ that verifies $Trace_{q}(\beta) = \alpha$ because, for such $\beta$,  the holomorphic $p-$form $\beta - \hat{q}^{*}(\frac 1kTrace_{q}(\beta))$ belongs to the kernel of $Trace_{q}$.
\hfill{$\square$}

\begin{lemma}\label{pull-back.5}
Consider a cartesian diagram of complex spaces
$$
\xymatrix{
\tilde{Z} \ar[d]_{\pi} \ar[r]^{\ \ \tilde{f}} & T \ar[d]^{q}\\ 
Z \ar[r]^{f} & S
} 
$$
where $q$ is surjective, proper, finite and geometrically flat and $Z$ is irreducible. Suppose also that we have an analytic cycle $G$ whose support is $\tilde{Z}$ and that $q$ admits a pull-back morphism $\hat{q}^* : \omega_S\rightarrow q_*\Omega_T$. Then $\frac 1{\deg \pi}Trace_{\pi}(\tilde{f}^*(\hat{q}^*(\alpha))) = \frac 1{\deg G}Trace_{G/Z}(\tilde{f}^*(\hat{q}^*(\alpha)))$ for all $\alpha\in\omega_M^{\bullet}(M)$.
\end{lemma}

\parag{Proof} If $f(Z)$ is not contained in the singular part of $S$ the result is obvious. If $f(Z)$ is contained in the singular part of $S$ we use the  same kind of arguments as in the proofs of lemmas \ref{pull-back.3} and \ref{pull-back.4}.
\hfill{$\blacksquare$}

\parag{Remark}
The lemma \ref{pull-back.5} tells us that the $\mathcal{O}_S-$morphism $\hat{f}_G^* : \omega_S^{\bullet}\rightarrow f_*\omega_Z^{\bullet}$ is independent of the analytic cycle $G$. (See lemma \ref{pull-back.-1} for the notation.)
\hfill{$\square$}

\begin{lemma}\label{pull-back.0}
Consider  the commutative diagram of complex spaces
$$
\xymatrix{
Z\times_{S} {T_2}\ar@/_3pc/[dd]_\pi \ar[d]_{\pi_2} \ar[r]^{\ \tilde{\tilde{f}}} & {T_2} \ar@/^3pc/[dd]^q\ar[d]^{q_2} \\ Z \times_{S} {T_1} \ar[r]^{\ \tilde{f}} \ar[d]_{\pi_{1}} & T_1 \ar[d]^{q_{1}}\\ 
Z \ar[r]^{f} & {S}} 
$$
where $S$ and $Z$ are irreducible, and $q_1 : T_1\rightarrow S$, $q_2 : T_2\rightarrow T_1$ are surjective, proper, finite and geometrically flat maps of degrees $k_1$, $k_2$ (resp.). Denote respectively $G_1$, $G_2$ and $G$ the graph cycles of  $q_1$, $q_2$ and $q$.  Then  for  any open subset $V$ of $S$ and any $\gamma\in\Omega_{T_1}^{\bullet}(q_1^{-1}(V))$ we have
\begin{equation}\label{eq.0}
\frac 1{k_1k_2}Trace_{f^*G/Z}(\tilde{\tilde{f}}^*(q_2^*\gamma)) \ = \ \frac 1{k_1}Trace_{f^*G_1/Z}(\tilde{f}^*(\gamma)).
\end{equation}
\end{lemma}

\parag{Proof}  We show first that for any open subset $U$ of $Z$ and any $\beta$ in $\Omega_{Z\times_ST_2}^{\bullet}(\pi^{-1}(U))$ one has
\begin{equation}\label{eq.1}
Trace_{f^*G/Z}(\beta) = Trace_{f^*G_1/Z}(Trace_{\tilde{f}^*G_2/Z \times_{S} T_1}(\beta)).
\end{equation}
To do so we first choose a nowhere dense analytic subset $\Sigma$ of $Z$ such that $S(Z)\subseteq \Sigma$, $\pi_1$ is a covering map over $Z\setminus \Sigma$ and $\pi_2$ is a covering map over $Z\times_ST_1\setminus \pi_1^{-1}(\Sigma)$. Since the problem is local and the $\mathcal{O}_Z-$module is $\omega_Z^{\bullet}$ is torsion free it is enough to prove the statement for all $U$ belonging to a basis of open sets in $Z\setminus \Sigma$. Now every point in $Z\setminus \Sigma$ has a basis of open connected neighbourhoods $U$ that are  $\pi_1-$trivializing, $\pi_1^{-1}(U) = V_1\cup\cdots\cup V_r$, and such that  each $V_j$ is $\pi_2-$trivializing, $\pi_1^{-1}(V_j) = W_{j1}\cup\cdots\cup W_{jl_j}$. The graph cycle $G_1$ endows each $V_j$ with a multiplicity $k_1^j$ such that $k_1^1 + \cdots + k_1^r = k_1$, and the graph cycle $G_2$ endows each $W_{jl}$ with a multiplicity $k_2^{jl}$ such that $k_2^{j1} + \cdots + k_2^{jl_j} = k_1^j$.

Then we obtain (\ref{eq.1}) by straightforward calculations.

\smallskip 
Now to prove (\ref{eq.0}) we take an open subset $V$ of $S$ and $\gamma$ in $\Omega_{T_1}^{\bullet}(q_1^{-1}(V))$. Then from (\ref{eq.1}) we deduce
\begin{align*}
&\frac 1{k_1k_2}Trace_{f^*G/Z}(\tilde{\tilde{f}}^*(q_2^*\gamma)) = \frac 1{k_1k_2}Trace_{f^*G_1/Z}(Trace_{\tilde{f}^*G_2/Z \times_{S} T_1}(\tilde{\tilde{f}}^*(q_2^*\gamma))) \\
&= \frac 1{k_1}Trace_{f^*G_1/Z}\left(\frac 1{k_2}Trace_{\tilde{f}^*G_2/Z \times_{S} T_1}(\pi_2^*(\tilde{f}^*\gamma))\right) = 
\frac 1{k_1}Trace_{f^*G_1/Z}(\tilde{f}^*\gamma)).
\end{align*}
\hfill $\blacksquare$

\begin{lemma}\label{pull-back.1}
Consider the diagram of complex spaces and holomorphic maps
$$
\xymatrix{
{} & T_1 \ar[rd]^{q_1}&&T_2 \ar[ld]_{q_2}\\ 
Z \ar[rr]^{f} && S& } 
$$
where $S$ and $Z$ are irreducible, $q_1$ and  $q_2$ are proper, finite and geometrically flat of degrees $k_1$, $k_2$. For $j=1,2$ we suppose that $q_j$ admits a pull-back morphism  $\hat{q}_j : \omega_S^{\bullet}\rightarrow(q_j)_*\Omega_{T_j}^{\bullet}$, we let $G_j$ denote the graph cycle of $q_j$ and $f_j : Z\times_ST_j\rightarrow T_j$ denote the natural map.
Then for every open subset $U$ of $S$ and every $\alpha$ in $\omega_S^{\bullet}(U)$ we have
\begin{equation*}
\frac 1{k_1}Trace_{f^*G_1/Z}(f_1^*(\hat{q}_1^*\alpha)) \ = \ \frac 1{k_2}Trace_{f^*G_2/Z}(f_2^*(\hat{q}_2^*\alpha)).
\end{equation*}
\end{lemma}
\parag{Proof} Put  $T := T_1\times_ST_2$ and apply lemma \ref{pull-back.0} to the diagram
$$
\xymatrix{
Z\times_S {T} \ar[d] \ar[r] & T \ar[d] \\ 
Z \times_{S} {T_j} \ar[r] \ar[d] & T_j \ar[d]^{q_{j}}\\ 
Z \ar[r]^{f} & S} 
$$
for \ $j = 1,2$.
\hfill $\blacksquare$

\begin{thm}\label{pull-back.2}
Let $f : Z\rightarrow M$ be a holomorphic map from a pure dimensional reduced complex space to a nearly smooth complex space. Then there exists a unique graded pull-back morphism 
$$
\hat{f}^* : \omega_M^{\bullet}\longrightarrow f_*\omega_Z^{\bullet}
$$
having the following properties.
\begin{enumerate}[i)]
\item
The diagram \ 
$
\xymatrix{
\Omega_{M}^{\bullet}  \ar[r]^{f^*}\ar[d]_{h_M^{\bullet}}& \ f_{*}\Omega_{Z}^{\bullet}\ar[d]^{h_Z^{\bullet}}\\ 
\omega_{M}^{\bullet}  \ar[r]^{\hat{f}^*} & \ f_{*}\omega_{Z}^{\bullet}} 
$
\ is commutative.
\item
Let $q : \tilde{M}'\rightarrow M'$ be a local model  of degree $k$ on $M$ and let $G$ denote its graph cycle. Then, for every $\alpha\in\omega_{M'}^{\bullet}(M')$
one has the identity 
$$
\hat{f}^*(\alpha) = \frac 1kTrace_{f^*G/Z}(\tilde{f}^*(\hat{q}^*\alpha))
$$ 
where  $\tilde{f} : f^{-1}(M')\times_{M'}\tilde{M}'\rightarrow \tilde{M}'$ is the natural map.
\item
The morphism $\hat{f}^*$ is a morphism of graded differential modules over the graded differential algebra $(\Omega_M^{\bullet},d^{\bullet})$.
\end{enumerate}
\end{thm}

\parag{Proof} We cover $M$ by local models and define a pull-back morphism  for each local model by the formula in  ii). Then by lemma \ref{pull-back.1} these pull-back morphisms glue together to the global morphism $\hat{f}^*$. The morphism verifies property  ii) by construction and its uniqueness is again due to lemma \ref{pull-back.1}. 

\smallskip
To show that $\hat{f}^*$ satisfies property i) we may suppose that we are in a local model;  that corresponds to the following diagram
$$
\xymatrix{
Z \times_{M}{\tilde{M}} \ar[r]^{\quad \tilde{f}} \ar[d]_{\tilde{q}} & \tilde{M} \ar[d]^{q}\\ Z \ar[r]^{f} & {M}} 
$$
Let $G$ be the graph cycle of $q$ and put $k := \deg q$. Then, for every $\alpha\in\Omega_{M}^{\bullet}(M)$, we get
\begin{eqnarray*}
\hat{f}^*(h_M^{\bullet}(\alpha)) 
&=& \frac 1kTrace_{f^*G/Z}(\tilde{f}^*\hat{q}^*(h_M^{\bullet}(\alpha)) 
= \frac 1kTrace_{f^*G/Z}(\tilde{f}^*q^*(\alpha))\\
&=& \frac 1kTrace_{f^*G/Z}(\tilde{q}^*f^*(\alpha)) \ 
= \ h_Z^{\bullet}(f^*(\alpha)).
\end{eqnarray*}
From ii) we immediately get iii).
\hfill{$\blacksquare$}

\parag{Remarks}
\begin{enumerate}[i)]
\item
If $Z$ is a locally closed analytic subset in  $M$ then the theorem tells us that there is a natural restriction morphism $\omega_M^{\bullet}\rightarrow \omega_Z^{\bullet}$. In particular if $Z$ is smooth there is a natural restriction morphism $\omega_M^{\bullet}\rightarrow \Omega_Z^{\bullet}$.
\item
From theorem \ref{pull-back.2} we see that if $M$ is a nearly smooth complex space and $\pi : M_{1}\to M$ is a desingularization of $M$, we have the identity of graded $\mathcal{O}_M-$modules
$$ 
\pi_{*}(\Omega_{M_{1}}^{\bullet}) = \omega_{M}^{\bullet}.
$$
As a nearly smooth complex space is normal and locally  $\mathbb{Q}-$Gorentein, it is Kawamata-log-terminal.
\end{enumerate}

The functoriality of the pull-back morphism $\hat{f}^{*}$ given by the next proposition will be useful.

\begin{prop}\label{funct. pull-back}
Let $f : Z \to M$ and $g : M \to N$ be holomorphic maps where $Z$ is a complex space and $M$ and $N$ are nearly smooth complex spaces. Then, for every open subset $U$ of $N$ and every $\alpha$ in $\omega_N^{\bullet}(U)$, we have $\widehat{g\circ f}^*(\alpha) = \hat{f}^*(\hat{g}^*(\alpha))$.
\end{prop}

\parag{Proof} As the assertion is local, we may assume that we have local models for $M$ and $N$, $\tilde{M} \rightarrow M$ and $p : \tilde{N} \to N$, and that $\alpha$ is a global section of \ $\omega_N^{\bullet}$. We may also suppose that $Z$ is irreducible, even smooth.
Consider the commutative diagram
$$
\xymatrix{
Z\times_{M} \tilde{M}\times_{N}\tilde{N} \ar[d] \ar[r]^{\ \ \tilde{\tilde{f}}} & \tilde{M}\times_{N} \tilde{N} \ar[d]^{q_2}& \quad  \\ 
Z \times_{N} \tilde{N} \ar[r]^{\tilde{f}} \ar[d] & M \times_{N}\tilde{N} \ar[d]^{q_1} \ar[r]^{\quad\tilde{g}} & \tilde{N} \ar[d]^{p}\\ 
Z \ar[r]^{f} & M \ar[r]^{g} & N
} 
$$
and let  $G_1$ and $G$ denote the graph cycles of $q_1$ and $q_1\circ q_2$ (resp.).  From lemma \ref{pull-back.5} we then have
$$
\widehat{g\circ f}^{*}(\alpha) \ = \ \frac 1{\deg q_1}Trace_{f^*G_1/Z}(\tilde{f}^*(\tilde{g}^*( \hat{p}^*(\alpha)))) \ = 
\ \frac 1{\deg q_1}Trace_{f^*G_1/Z}(\tilde{f}^*\gamma)
$$
where \ $\gamma := \tilde{g}^*( \hat{p}^*(\alpha)))$. On the other hand we also get from lemma \ref{pull-back.5} that 
$$
\hat{f}^*(\hat{g}^*(\alpha)) = \ \hat{f}^*\left(\frac 1{\deg q_1}Trace_{q_1}(\gamma)\right). 
$$
Since the map $\tilde{M}\times_{N} \tilde{N}\rightarrow M$ factorizes through the map $\tilde{M} \rightarrow M$ we obtain by lemma \ref{pull-back.0} and the remark made after lemma \ref{pull-back.4}:
$$
\hat{f}^*\left(\frac 1{\deg q_1}Trace_{q_1}(\gamma)\right) = \frac 1{\deg q_2\cdot\deg q_1}Trace_{f^*G/Z}(\tilde{\tilde{f}}^*(q_2^*\gamma))
$$
Then again from lemma \ref{pull-back.0} we finally get 
$$
\frac 1{\deg q_2\cdot\deg q_1}Trace_{f^*G/Z}(\tilde{\tilde{f}}^*(q_2^*\gamma)) \ = \ \frac 1{\deg q_1}Trace_{f^*G_1}(\tilde{f}^*\gamma)
$$
\hfill{$\blacksquare$}\\

\section{ Fundamental classes}

\subsection{Definition and first properties.}

\parag{Notation and terminology} 
\begin{itemize}
\item 
In this section a $n-$cycle in $M$ means a locally finite sum of irreducible $n-$dimensional analytic subsets with positive {\em integers} as coefficients.
\item
Recall that for a geometrically flat map $f : M \to N$ between complex spaces and a cycle $Y$ in $N$ there exists a natural pull-back cycle of $Y$ by $f$ in $M$, denoted $f^*Y$. It is additive in $Y$ and, for a given irreducible  $Y$, it is defined as the graph cycle  in $N \times_{N} M \simeq M$ of the analytic family of fibres of $f$ parametrized by $Y$. It is proved in  [B-M 2]\, ch.VI \, that in this situation the pull-back of an analytic family of cycles in $N$ is an analytic family of cycles in $M$. 
\item 
For an analytic subset $X$ in a complex space $M$ we shall denote $\underline{H}_{Xmod}^{s}(\omega_{M}^{r}) $ the submodule of the \ $\mathcal{O}_M-$module $\underline{H}_{X}^{s}(\omega_{M}^{r})$ that is annihilated (locally) by some power of the reduced ideal sheaf $\mathcal{I}_{X}$ of $X$ in $M$. It is called the {\bf moderate cohomology sheaf} with support in $X$ and coefficients in $\omega_{M}^{r}$.
\end{itemize}

In the case where  $M$ is smooth we have $\Omega_{M}^{r} = \omega_{M}^{r}$ and these moderate cohomology sheaves are the cohomology sheaves of the $\bar\partial-$complex of $(r, \bullet)-$currents with support in $X$. If we suppose moreover that \ $M$ \ is connected and that \ $X$ \ is a cycle of  codimension \ $p$ \ in \ $M$, then its integration current \ $[X]$ \ defines a global section of \ $\underline{H}_{Xmod}^{p}(\omega_{M}^{p})$, called the {\bf fundamental class} of \ $X$ \ in \ $M$ and denoted by \ $c_X^M$ \  (see [B-M 2]). 

\medskip
Now suppose that we have a local model $q : \tilde{M} \to M$ and a cycle $Y$ of codimension \ $p$ \ in $M$. Put \ $X := q^{*}(Y)$. The functor $q_{*}$ being exact we have a natural isomorphism
$$
q_{*}\left(\underline{H}_{|X|}^{p}(\Omega_{\tilde{M}}^p)\right) \longrightarrow \underline{H}_{|Y|}^{p}\left(q_{*}(\Omega_{\tilde{M}}^p)\right)
$$
so we can identify the two sheaves. By composing this isomorphism with the morphism  
$$
\underline{H}_{|Y|}^{p}\left(q_{*}(\Omega_{\tilde{M}}^p)\right)\longrightarrow \underline{H}_{|Y|}^{p}\left(\omega_{M}^p\right)
$$
induced by the trace map $Trace_{q} : q_{*}(\Omega_{\tilde{M}}^{p}) \to \omega_{M}^{p}$ we get an \ $\mathcal{O}_M-$linear map (which we will also denote by $Trace_{q}$)
\begin{equation*}
Trace_{q} : \ q_{*}\left(\underline{H}_{|X|}^{p}(\Omega_{\tilde{M}}^p)\right) \longrightarrow \underline{H}_{|Y|}^{p}\left(\omega_{M}^p\right).
\end{equation*}

\parag{Notation}
Let $Y$ be a $n-$cycle in a nearly smooth complex space $M$ of pure dimension $n+p$ and let $q : \tilde{M}' \rightarrow M'$ be a local model of degree $k$ for $M$. Then  the section on $M'$ of the  sheaf \ $\underline{H}_{|Y|}^{p}\left(\omega_{M}^p\right)$ given by
$$ 
c_{Y\cap M'}^{M'}(q) := \frac 1k.Trace_{q}(c_{X}^{\tilde M}), 
$$
where $X := q^*(Y\cap M')$, will be called the (analytic)  {\bf fundamental class of $Y$ with respect to $q$} (or in the local model).
Notice that the section  $c_{Y\cap M'}^{M'}(q)$  is annihilated by the reduced ideal of $|Y|$   in $\mathcal{O}_{M'}$ because the reduced ideal  of \ $|X|$ \ in $\mathcal{O}_{\tilde M'}$ annihilates $c_{X}^{\tilde M}$. In particular we see that $c_{Y\cap M'}^{M'}(q)$ is a section of the moderate cohomology sheaf.

\begin{thm}\label{class. fond.}
Let $Y$ be a cycle of codimension $p$ in a nearly smooth complex space $M$. Then there exists a unique global section $c_{Y}^{M}$ of the sheaf $\underline{H}_{|Y|mod}^{p}(\omega_{M}^{p})$ such that for every local  model $q :  \tilde{M'} \rightarrow M'$ for $M$ the restriction of $c_{Y}^{M}$ to $M'$ coincides with $c_{Y\cap M'}^{M'}(q)$. This section is $d-$closed and is annihilated by the reduced ideal of $\vert Y\vert$.
\end{thm}

\begin{defn}
In the situation of theorem \ref{class. fond.} the  section $c_{Y}^{M}$ will be called the (analytic) {\bf\em fundamental class} of $Y$ in $M$
\end{defn}

\parag{Proof of theorem \ref{class. fond.}} Since $M$ can be covered by local models it is sufficient to show that different local models on the same open subset $V$ of $M$ will give the same fundamental class for $Y$ on $V$. To do so we will assume that we have a local model $q :  \tilde{M} \to M$ and prove that $c_{Y}^{M}(q)$ does not depend on $q$.
Since the problem is local it is enough to consider the case when $Y$ is an irreducible analytic subset of $M$. If $Y$ is not contained in the singular part $S(M)$ of $M$, then  the restriction of $c^M_Y(q)$ to $M\setminus S(M)$ is determined by the integration current of $Y\setminus S(M)$  and consequently independent of $q$. Hence $c^M_Y(q)$ is independent of $q$, thanks to proposition \ref{N-Smooth 2} v).\\
Now consider the case where $Y$ is contained in the singular part of $M$. Again we know from proposition \ref{N-Smooth 2}  v) that the sheaf $\underline{H}_{Y}^{p}(\omega_{M}^{p})$ has no non zero section whose support is contained in a nowhere dense analytic subset of  $Y$, so it is enough to establish the result in the case when $Y$ is a smooth submanifold of the singular part of $M$. Let $l$ denote the lowest multiplicity  in $M$ of a point $y$ in $Y$. Then the set of points in $Y$ that have multiplicity $l$ in $M$ contains a dense Zariski open subset in $Y$. Let  $y_{0}$ be one of these points. Then by shrinking $M$ to a smaller open neighbourhood of $y_0$ we may assume that we have a map $\pi : M \rightarrow U$ (local parametrization) onto an open polydisc $U$ in $\C^{m}$ with the following properties: 
\begin{enumerate}[(i)]
\item 
The map $\pi : M \to U$ is a surjective proper and finite of degree $l$.
\item 
The image of $Y$ is a connected sub-manifold $Y_{0}$ in $U$.
\item 
The pull-back cycle $\pi^{*}(Y_{0})$ is $l.Y$.
\end{enumerate}
Note that these conditions imply that the pull-back cycle $X :=  q^{*}(Y) $ satisfies $\pi_{*}(X) = k.l.Y_{0}$ as $q_{*}(X) = k.Y$ and $\pi_*(Y) = l.Y_{0}$.

These conditions may clearly be satisfied near each point of  on an open dense set in $Y$ on which a section of the sheaf $\underline{H}_{X}^{p}(\omega_{M}^{p})$ is determined\footnote{This open set is constructed by removing finitely many times a nowhere dense analytic subset  in the previous open dense subset. A section of the sheaf $\underline{H}_{X}^{p}(\omega_{M}^{p})$ is determined by its restriction to an open set of this kind.}. 
So now we have the following situation:
\begin{itemize}
\item 
Two holomorphic surjective proper and finite maps:
$$ q : \tilde{M} \to M \quad {\rm and} \quad \pi : M \to U$$
 of degrees respectively equal to $k$ and $l$. 
\item 
We have traces $Trace_{q}: q_{*}(\Omega_{\tilde{M}}^{p}) \to \omega_{M}^{p}$, \ $Trace_{\pi}: \pi_{*}(\omega_{M}^{p}) \to \Omega_{U}^{p}$ \ and  \ $Trace_{\pi\circ q}: (\pi\circ q)_{*}(\Omega_{\tilde{M}}^{p}) \to \Omega_{U}^{p}$.
\end{itemize}

As $(\pi\circ q)^{*}(Y_{0}) = l.X$  the compatibility of the pull-back of cycles with the pull-back of fundamental classes for a geometrically flat map between complex manifolds (see [B-M 2] chapter VI) gives
\begin{equation*}
Trace_{q}(c_{l.X}^{\tilde{M}}) = Trace_{q}((\pi\circ q)^{*}(c_{Y_{0}}^{U})) = Trace_{q}(q^{*}({\pi}^{*}(c_{Y_{0}}^{U})))  = k.{\pi}^{*}(c_{Y_{0}}^{U})
\end{equation*}
because $Trace_{q}\circ q^{*} = k.id_{\omega_{M}^{\bullet}}$. So we obtain
\begin{equation}\label{eq.2}
 \frac{1}{k}.Trace_{q}(c_{X}^{\tilde{M}}) = \frac{1}{l}.{\pi}^{*}(c_{Y_{0}}^{U}) 
 \end{equation}
and this completes the proof as the right-hand side of the formula (\ref{eq.2}) does not depend on  $q$ but only of the choice of the local parametrization $\pi$ of $M$.
$\hfill\blacksquare$\\

We remark that the left-hand side of formula (\ref{eq.2}) does not depend on the choice of the local parametrisation $\pi : M \to U$, so it gives an alternative computation of the fundamental class of the cycle $X$ near a smooth generic point of  $|X|$ when $\vert X\vert$ is contained in the singular set of $M$. This is a very special case of the theorem \ref{pull-back fund. class} below.
 
\parag{Some properties of the fundamental class}
Let $X$ be a cycle of codimension $p$ in a nearly smooth complex space $M$.
\begin{itemize}
\item 
As in the smooth case the fundamental class $c_{X}^{M} \in \Gamma(M,\underline{H}_{ Xmod}^{p}(\omega_{M}^{p}))$  is  $d-$closed, annihilated by $\mathcal{I}_{\vert X \vert}$  and consequently annihilated by exterior product by $df$ for any $f \in \mathcal{I}_{|X|}$.
\item 
The exactness in positive degrees of the de Rham complex $(\omega_{M}^{\bullet}, d^{\bullet})$ and the vanishing on $M$  of the sheaves $\underline{H}_{|X|}^{j}(\omega_{M}^{\bullet})$ for all $j \not = p$  show that we have a canonical isomorphism of sheaves 
$$
top_M : H^{p}\big(\underline{H}_{|X|}^{p}(\omega_{M}^{\bullet}), d^{\bullet}\big) \to \underline{H}_{|X|}^{2p}(\underline{\C}).
$$
The image of the fundamental class $[c_{X}^{M}]$ by this morphism will be called {\bf the topological fundamental class of $X$ in $M$}. 

Now let $q : \tilde{M'} \rightarrow M'$ be a local model of degree $k$ on $M$. Then  the following diagram is commutative
$$
\xymatrix{ 
q_*\left(H^{p}\left(\underline{H}_{|q^*X|}^{p}(\Omega_{\tilde M'}^{\bullet}), d^{\bullet}\right)\right) \ar[dd]_{Trace_q} \ar[rr]^{\qquad top_{\tilde M'} }&& q_*\left(\underline{H}_{|q^*X|}^{2p}(\underline{\C}_{\tilde M'})\right)\ar[dd]^{Trace_q^{top}}\\
\\ 
H^{p}\big(\underline{H}_{|X|}^{p}(\omega_{M'}^{\bullet}), d^{\bullet}\big)  \ar[rr]^{\qquad top_M} && \underline{H}_{|X|}^{2p}(\underline{\C}_{M'}) }
$$
where the morphism $Trace_q^{top}$ is induced by the morphism $\underline{\C}_{\tilde M'}\rightarrow\underline{\C}_{M'}$ that, for every open set $V$ in $M'$, maps a locally constant function $f : q^{-1}(V)\rightarrow\C$ to the locally constant function 
$$
V\longrightarrow \C,\qquad x\mapsto \sum_{y\in q^*(x)}f(y).
$$
Hence the restriction to $M'$ of the topological fundamental class $top_M([c^M_X])$ is given by the formula
\begin{equation}\label{eq.A}
top_M([c^M_X])_{|_{M'}} \ = \ \frac 1kTrace_q^{top}\left(top_{\tilde M}([c^{\tilde M}_{q^*X}])\right).
\end{equation}
In a complex manifold the topological fundamental class of an analytic cycle $Z$ is a global section of $\underline{H}_{|Z|}^{2p}(\underline{\mathbb{C}})$ that takes its values in the image sheaf of the canonical morphism $\underline{H}_{|Z|}^{2p}(\underline{\mathbb{Z}})\rightarrow \underline{H}_{|Z|}^{2p}(\underline{\mathbb{C}})$. From formula (\ref{eq.A}) we see that this need not be true when $M$ is singular and we can only affirm that this section takes its values in the image sheaf of the morphism  $\underline{H}_{|Z|}^{2p}(\underline{\mathbb{Q}})\rightarrow \underline{H}_{|Z|}^{2p}(\underline{\mathbb{C}})$.
This fact should be seen in the context of section 3 where we use rational cycles to develope a reasonable intersection theory for nearly smooth spaces.
\hfill{$\square$}
\end{itemize}
  
\begin{lemma}\label{caract.}
Let $M$ be a nearly smooth complex space and let $X$ be a cycle of codimension $p$ in $M$. The fundamental class $c_{X}^{M}$ of $M$ is the unique section of the sheaf $\underline{H}_{\vert X\vert}^{p}(\omega_{M}^{p})$ with the following properties:
\begin{enumerate}[(i)]
\item 
$dc_{X}^{M} = 0$ as a section of the sheaf $\underline{H}_{\vert X\vert}^{p}(\omega_{M}^{p+1})$;
\item 
$\mathcal{I}_{\vert X\vert}.c_{X}^{M} = 0$ where $\mathcal{I}_{\vert X\vert}$ is the reduced ideal sheaf of the support $\vert X\vert$ of $X$;
\item 
The image of $c_{X}^{M}$ in $\underline{H}_{\vert X\vert}^{2p}(\C)$ is the topological fundamental class of $X$.
\end{enumerate}
\end{lemma}
  
\parag{Proof} As this statement is local, it is an immediate consequence of the analogous statement in the smooth case.
$\hfill \blacksquare$
  
\subsection{Relative fundamental classes.}

Let $S$ and $M$ be two complex spaces and let $p_{M}: S\times M \to M$ denote the natural projection. We put
$$
\omega_{S\times M/S}^{\bullet} := p_{M}^{*}(\omega_{M}^{\bullet}).
$$ 
For each $s \in S$ let $j_{s} : M \to S \times M$ denote the injection given by $x \mapsto (s, x)$. Then there exists a natural pull-back morphism $j_{s}^{*} : \omega^{\bullet}_{S\times M/S} \to (j_{s})_{*}(\omega_{M}^{\bullet})$.

Now suppose we have a local model $q : \tilde{M}\rightarrow M$ and consider the commutative diagram
$$
\xymatrix{ S \times \tilde{M} \ar[d]_{\id_{S}\times q} \ar[r]^{\tilde{p}} & \tilde{M} \ar[d]^{q}\\ S \times M \ar[r]^{p} & M }.
$$
The trace map $Trace_q$ induces an $\mathcal{O}_{S\times M}-$linear morphism 
$$
p^*(Trace_q) : p^{*}(q_{*}(\Omega_{\tilde{M}}^{\bullet})) \rightarrow p^{*}(\omega_{M}^{\bullet}).
$$
and by definition we have 
$$
(\id_{S}\times q)_{*}(\Omega_{S \times \tilde{M}/S}^{\bullet}) = (\id_{S}\times q)_{*}(\tilde{p}^{*}(\Omega_{\tilde{M}}^{\bullet}))
 $$
and $p^{*}(\omega_{M}^{\bullet}) = \omega_{S\times M/S}^{\bullet}$.
Since $q$ is a proper finite map the natural morphism
$$
(\id_{S}\times q)_{*}(\tilde{p}^{*}(\Omega_{\tilde{M}}^{\bullet})) \ \longrightarrow \  p^{*}(q_{*}(\Omega_{\tilde{M}}^{\bullet})) 
$$
is an isomorphism\footnote{This reduces to showing that the natural map 
$$
(\mathcal{O}_{S,s}\hat{\otimes}_{\C}\mathcal{O}_{M,q(x)})\otimes_{\mathcal{O}_{M,q(x)}}\mathcal{O}_{\tilde{M},x} \longrightarrow \mathcal{O}_{S,s}\hat{\otimes}_{\C}\mathcal{O}_{\tilde{M},x}
$$
is an isomorphism for every $(s,x)$ in $S\times \tilde{M}$, which is straightforward. Here $\hat{\otimes}_{\C}$ denotes the completed tensor product.} and by composing it with $p^*(Trace_q)$ we get the {\bf relative trace map}
$$
Trace_{q/S} : (\id_{S}\times q)_{*}(\Omega_{S \times \tilde{M}/S}^{\bullet}) \longrightarrow \omega_{S\times M/S}^{\bullet}.
$$
In the same way the pull-back morphism $\hat{q} :  \omega_{M}^{\bullet} \rightarrow q_{*}(\Omega_{\tilde{M}}^{\bullet})$ induces a {\bf relative pull-back morphism}
$$
\hat{q}_{/S} :  \omega_{S\times M/S}^{\bullet} \ \longrightarrow \ (\id_{S}\times q)_{*}(\Omega_{S \times \tilde{M}/S}^{\bullet}).
$$
It is an $\mathcal{O}_{S\times M}-$linear section of $Trace_{q/S}$ and consequently induces an isomorphism of  $\omega_{S\times M/S}^{\bullet}$ onto a direct factor of the $\mathcal{O}_{S\times M}-$module $(\id_{S}\times q)_{*}(\Omega_{S \times \tilde{M}/S}^{\bullet})$.

\smallskip
Let $U$ be an open subset of $M$ and, for each $s$ in $S$,  let $\tilde{j}_s : \tilde{M}\rightarrow S\times\tilde{M}$ denote the injection $x \mapsto (s, x)$. Then for every $\alpha$ in $\Gamma(S\times U,\omega_{S\times M/S}^{\bullet})$ and every $\beta$ in $\Gamma((\id_S\times q)^{-1}(S\times U),\Omega_{S\times \tilde{M}/S}^{\bullet})$ we have
\begin{equation}\label{eq.3}
\hat{q}^*(j_s^*(\alpha)) = \tilde{j}_s^*(\hat{q}_{/S}^*(\alpha))\qquad\text{and}\qquad j_s^*(Trace_{q/S}(\beta)) = Trace_q(\tilde{j}_s^*(\beta))
\end{equation}

In the sequel we will use the same symbols for the morphisms induced by the above morphisms on cohomology sheaves.

\begin{lemma}\label{unique}
Let $S$ be a complex space,  $M$ a nearly smooth complex space and  $Z$  an analytic subset in $S \times M$. Suppose $Z_s := p_M(Z\cap (\{s\}\times M))$ is of pure codimension $p$ for every $s$ in $S$. 
Let $c$ be a global section of the sheaf  $\underline{H}_{Z}^{p}(\omega_{S\times M/S}^{\bullet})$ such that for each $s \in S$ the pull-back section $(j_{s})^{*}(c)$ of the sheaf $\underline{H}_{Z_{s}}^{p}(\omega_{M}^{\bullet})$ is identically zero. Then the section $c$ is identically zero.
\end{lemma}

\parag{Proof} 
As the assertion is local on $S \times M$ we my assume that we have a local model $q :\tilde{M} \rightarrow M$. We know that the lemma is valid in the smooth case (see [B-M.2]) so from (\ref{eq.3}) we obtain $\hat{q}_{/S}^*(c) = 0$ and consequently $c = 0$ since $\hat{q}_{/S}^*$ induces  an isomorphism onto its image.
$\hfill \blacksquare$\\

\begin{defn}\label{class. fond.rel}
Let $S$ and $M$ be two complex spaces and assume that $M$ is nearly smooth. Let $(X_{s})_{s \in S}$ be a family of cycles of (pure) codimension $p$ in $M$ such that its set theoretic graph
$$ 
|\mathcal{X}| := \{ (s,t) \in S \times M \ / \  t \in |X_{s}| \} 
$$
is a analytic subset in $S \times M$. We shall say that a section $c $ of the sheaf $\underline{H}_{\vert \mathcal{X}\vert}^{p}(\omega_{S\times M/S}^{p})$ is a {\bf\em relative fundamental class} for the family $(X_{s})_{s \in S}$ if for each $s \in S$ the global section  $(j_{s})^{*}(c)$ of $\underline{H}_{\vert X_{s}\vert}^{p}(\omega_{M}^{p})$ is the fundamental class of the cycle $X_{s}$.
\end{defn}

\parag{Remark} If such a relative  fundamental class exists for the family $(X_{s})_{s \in S}$ it is unique, thanks to the lemma \ref{unique}.
 
\begin{thm}\label{class. fond. rel.}
Let $S$ and $M$ be two reduced complex spaces and assume that $M$ is nearly smooth. Let $(X_{s})_{s \in S}$ be a family of cycles of (pure) codimension $p$ in $M$ such that its set theoretic graph
$$ 
|\mathcal{X}| := \{ (s,t) \in S \times M \ / \  t \in \vert X_{s}\vert \} 
$$
is a analytic subset in $S \times M$. Then the existence of a relative fundamental class $c_{\mathcal{X}/S}^{M}$ for such a family of cycles in $M$  is equivalent to the analyticity of this family
\end{thm}

\parag{Proof} As the assertion is local, on $S \times M$ we may assume that we have a local model $q : \tilde{M} \to M$. The map $q$ is geometrically flat so, if the family $(X_{s})_{s \in S}$ is analytic in $M$, so is the family $(q^{*}(X_{s}))_{s \in S}$ in $\tilde{M}$ due to the pull-back theorem in [B-M 2] ch.VI . Conversely, if the family $(q^{*}(X_{s}))_{s \in S}$ is analytic in $\tilde{M}$ so is the family $(X_{s})_{s \in S}$ thanks to the direct image theorem (see [B-M 1] ch.IV) and the fact that $q_{*}(q^{*}(Y)) = k.Y$ for any cycle $Y$ in $M$. Also the analyticity of the set theoretic graph of the family $(X_{s})_{s \in S}$ or of the family $(q^{*}(X_{s}))_{s \in S}$ are equivalent. The pull-back and push-forward by $q$ give the equivalence of the existence of relative fundamental class for these families, completing the proof.
$\hfill \blacksquare$\\

The following theorem shows that, for a geometrically flat map  between two nearly smooth complex spaces, the pull-back of cycles is compatible with the pull-back of fundamental classes.

\begin{thm}\label{pull-back fund. class}
Let $f : M \rightarrow N$ be a geometrically flat map (see definition \ref{g-flat}) between two nearly smooth complex spaces and let  $Y$ be a cycle of codimension $p$ in $N$. Denote also 
$$
\hat{f}^{*} : \underline{H}_{|Y|}^{p}(\omega^{\bullet}_{N})\longrightarrow f_*\Big(\underline{H}_{|f^*Y|}^{p}(\omega^{\bullet}_{M})\Big)
$$
the morphism induced by the pull-back morphism $\hat{f}^{*} : \omega^{\bullet}_{N} \rightarrow f_{*}(\omega_{M}^{\bullet})$ (see theorem \ref{pull-back.2}). Then we have
$$ 
c_{f^*Y}^{M} =  \hat{f}^{*}(c_{Y}^{N}).
$$
\end{thm}

\parag{Proof} We know from [B-M 2] that the theorem is correct if both $M$ and $N$ are smooth so we are first going to consider the case where $N$ is smooth. The question being local on $M$  we may assume that we have a local model $q : \tilde{M} \to M$ of degree $k$ and since the map $f \circ q$ is geometrically flat we get
$$
(f\circ q)^*\left(c_Y^N\right) = c_{(f\circ q)^*Y}^{\tilde{M}} = c_{q^*(f^*Y)}^{\tilde{M}}.
$$
Due to proposition \ref{funct. pull-back} we have $(f\circ q)^* = \hat{q}^*\circ\hat{f}^*$ \ and hence $c_{q^*(f^*Y)}^{\tilde{M}} = \hat{q}^*\left(\hat{f}^*\left(c_Y^N\right)\right)$. By applying $\frac{1}{k}.Trace_{q}$ to both sides of this last identity we then obtain
$$
c_{f^*Y}^M = \frac{1}{k}.Trace_{q}\left(c_{q^*(f^*Y)}^{\tilde{M}}\right) = \hat{f}^*(c_Y^N)
$$
because $\frac{1}{k}.Trace_{q}\circ \hat{q}^{*} = \id_{\omega_M^{\bullet}}$.

In the general case, we may suppose $Y$ irreducible. If $Y$ is not contained in the singular part of $N$, the previous case is enough to conclude because the fundamental classes are determined by their restriction to an open subset meeting $Y$. If $Y$ is contained in the singular set of $N$ then, as the question is also local on $N$, we may assume that we have  a local parametrization $\pi : N \rightarrow U$,  such that $l.Y = \pi^{*}(Y_{0})$ where $Y_{0}$ is a smooth connected submanifold of the open polydisc $U$ and $l$ is the degree of $\pi$ (see the proof of the theorem \ref{class. fond.}). The map $\pi$ is geometrically flat (flat in fact, see the proposition \ref{N-Smooth 2} point ii)), so we may apply the previous case to the map $\pi\circ f$ to conclude because $c_{Y}^{N} = \frac{1}{l}.\hat{\pi}^{*}(c_{Y_{0}}^{U})$ and $\widehat{\pi\circ f}^{*} = \hat{f}^*\circ \hat{\pi}^{*}$ thanks to proposition \ref{funct. pull-back}.
$\hfill \blacksquare$

\begin{cor}\label{relative pull-back}
In the situation of the theorem \ref{pull-back fund. class}, if we have an analytic family $(Y_{s})_{s \in S}$ of  cycles of (pure) codimension $p$ in $N$ parametrized by a complex space $S$, the pull-back family $(f^{*}(Y_{s}))_{s \in S}$ is again an analytic family, and its relative fundamental class is the $S-$relative pull-back by $\id_{S}\times f$  of the relative fundamental class of the analytic  family $(Y_{s})_{s \in S}$.
\end{cor}

\parag{Proof} Combining the theorem \ref{pull-back fund. class} with the theorem \ref{class. fond. rel.} we obtain  the analyticity of the pull-back of the analytic family of cycles in $N$ 
\hfill$\blacksquare$

\section{Geometric intersection theory in a nearly smooth complex space}

\subsection{Preliminaries}

Our aim is now to establish a geometric intersection theory for analytic cycles in a nearly smooth complex space that is compatible with the smooth case. Some preliminaries are necessary and we begin by introducing the notion of a {\em rational cycle}.

\begin{defn}\label{rat. cycle}
Let $M$ be a complex space. A {\bf\em rational $n-$cycle}  $X$ in $M$ is a formal sum
$$ 
X := \sum_{j \in J} \  r_{j}.X_{j} 
$$
where the $r_{j}$ are positive rational numbers and  $(X_{j})_{j \in J}$ is the locally finite collection of irreducible components of a closed  analytic subset $ \vert X\vert$ of pure dimension $n$ in $M$. This subset $\vert X\vert$ will be called the { \bf\em support} of the cycle $X$. If the sum is empty, the corresponding cycle is by convention the  empty $n-$cycle in $M$ and its support is the empty set.
\end{defn}

The analytic cycles in the ordinary sense are the rational cycles whose coefficients are all integral. In this context we will call them  {\bf integral cycles}.\\

When we consider a family of rational cycles in a complex space $M$ the following condition will be important.

\begin{defn}\label{cond. D}
We shall say that a family $(X_{s})_{s \in S}$ of rational cycles in a complex space $M$ parametrized by a Hausdorff topological space $S$ satisfies the {\bf\em condition} $(D)$ on an open set $S' \times M'$ in $S \times M$ if there exists  an integer $d \geq 1$ such that for any $s \in S'$ the cycle $d.X_{s}\cap M'$ is integral. We shall say that the family $(X_{s})_{s \in S}$ satisfies the condition $(D)$ on $M$ if for any $(s, x) \in S \times M$ there exists open neighbourhoods $S'$ and $M'$ respectively  of $s$ in $S$ and of $x$ in $M$ such that the condition $(D)$ is satisfied on $S' \times M'$.
\end{defn}

\begin{defn}\label{anal. rat.}
Let $S$ be a complex space. We shall say that a family $(X_{s})_{s \in S}$ of rational cycles in a complex space $M$ parametrized by $S$ is {\bf\em  analytic} (resp. {\bf\em continuous)}  if it satisfies the condition $(D)$ and if, for any open set $S' \times M'$ and any integer $d$ such that $d.X_{s}\cap M'$ is an integral cycle for all $s \in S'$ the family $(d.X_{s}\cap M')_{s \in S'}$ of integral cycles is analytic (resp. continuous).
\end{defn}

\parag{Important remark}{\em Using the previous definition, it is immediate to extend all results obtained in the sections 1 and 2 to the case of rational cycles in nearly smooth complex spaces}.
\hfill{$\square$}\\

Let $X$ and $Y$ be two irreducible analytic subsets in a nearly smooth complex space $M$. Then by proposition \ref{N-Smooth 2} iii) we have
$$ 
\codim_{M}X \cap Y \leq \codim_{M}X + \codim_{M}Y
$$
and from this we see that the equality $\codim_{M}X \cap Y = \codim_{M}X + \codim_{M}Y$ implies that  $X\cap Y$ is of pure codimension.

\begin{defn}\label{coupe bien}
Let $M$ be a nearly smooth complex space and let $X$ and $Y$ be two rational cycles in $M$ of (pure) codimensions  $p$ and $p'$. We shall say that $X$ and $Y$  {\bf\em  intersect properly in $M$} when the intersection $|X|\cap|Y|$ is either empty or when it has (pure) codimension $p + p'$.
\end{defn} 

We notice that  if $X$ and $Y$ intersect properly each irreducible component of $|X|$ and each irreducible component of $|Y|$ intersect properly.

\begin{thm}\label{intersect. 1}
Let $M$ be a nearly smooth complex space and let $X$ and $Y$ be two integral cycles that  intersect properly  in $M$. Then there exist a unique rational cycle $X\cap_MY$ in $M$ such that for every local model $q : \tilde{M'}\rightarrow M'$  for $M$ we have
$$
(X\cap_MY)\cap M' = \frac{1}{\deg q}.q_{*}(q^{*}X \cap q^{*}Y)\cap M'.
$$
\end{thm}

\parag{Proof} Without loss of generality we may suppose that $X$ and $Y$ are irreducible analytic subsets of $M$. Since $M$ can be covered with local models it is enough to prove that for two local models $q_1 : M_1\rightarrow M'$ and $q_2 : M_2\rightarrow M'$ of degrees $k_1$, $k_2$ we have
$$
\frac{1}{k_1}.(q_1)_{*}(q_1^{*}X \cap q_1^{*}Y) = \frac{1}{k_2}.(q_2)_{*}(q_2^{*}X \cap q_2^{*}Y).
$$
We notice first that on the regular part $M'_{reg}$ of $M'$ the intersection cycle $X\cap_{M'_{reg}}Y$ is well defined and according to the {\em projection formula} for manifolds (see [B-M.2] chapter VII) we have
$$
\frac{1}{k_i}.(q_i)_{*}(q_i^{*}X \cap q_i^{*}Y) = \frac{1}{k_i}.(q_i)_{*}(q_i^{*}X) \cap Y = X\cap_{M'_{reg}}Y
$$
for $i = 1,2$ on $M'_{reg}$. Hence, if no irreducible component of $X\cap Y$ is contained in the singular part $S(M')$ the result holds, since in that case the cycle  $\frac{1}{k_i}.(q_i)_{*}(q_i^{*}X \cap q_i^{*}Y)$ is determined by its restriction to $M'_{reg}$.
\\
In the general case we can, by shrinking $M'$ if necessary, find two analytic families of integral cycles $(Z_s)_{s\in S}$  and $(W_t)_{t\in T}$ in $M_1$, parametrized by open neighbourhoods of the origins $S$ and $T$ in some numerical spaces having the following properties  (see [B-M.2] chapter VII).
\begin{itemize}
\item
$Z_0 = q_1^*X$ and  $W_0 = q_1^*Y$.
\item 
For all $(s,t)$ the cycles $Z_s$ and $W_t$ intersect properly.
\item
For all $(s,t)$ in an open dense subset of $S\times T$ no irreducible component of $|Z_s|\cap |W_t|$ is contained in $q_1^{-1}(S(M'))$.
\end{itemize}
Then $((q_1)_*Z_s)_{s\in S}$  and $((q_1)_*W_t)_{t\in T}$  are analytic families of cycles in $M'$ such that  $(q_1)_*Z_0 = k_1X$, $(q_1)_*W_0 = k_1Y$.  Then the rational cycles $X_s := \frac 1{k_1}(q_1)_*Z_s$ and $Y_t := \frac 1{k_1}(q_1)_*W_s$ form analytic families of rational cycles such that $X_0 = X$ and $Y_0 = Y$. It follows that the two families of cycles 
$$
\left(\frac{1}{k_1}.(q_1)_*(q_1^{*}(X_{s}) \cap q_1^{*}(Y_{s}))\right)_{(s,t)\in S\times T}\qquad\text{and}\qquad
\left(\frac 1{k_2}.(q_2)_*(q_2^{*}(X_{s}) \cap q_2^{*}(Y_{s}))\right)_{(s,t)\in S\times T}
$$
are both analytic and coincide on an open dense subset of $S\times T$. Consequently they are identical and in particular for $(s,t) = (0,0)$ we get
$$
\frac{1}{k_1}.(q_1)_{*}(q_1^{*}X \cap q_1^{*}Y) = \frac{1}{k_2}.(q_2)_{*}(q_2^{*}X \cap q_2^{*}Y).
$$
This completes the proof.
\hfill{$\blacksquare$}

\begin{defn}\label{intersect.2}
In the situation of theorem \ref{intersect. 1} the rational cycle $X\cap_MY$ is called the {\bf\em intersection cycle} of $X$ and $Y$ in $M$. 
\end{defn}

\parag{Remarks}
\begin{enumerate}[i)]
\item 
Theorem \ref{intersect. 1} provides us with a geometric intersection theory for rational analytic cycles on nearly smooth complex spaces. This theory is local, i.e.  for any open set $M'$ in a nearly smooth space $M$ and every pair $(X, Y)$ of cycles that intersect properly in $M$ we have
$$ 
(X \cap_{M} Y)\cap M' = (X \cap M')\cap_{M'}(Y \cap M').
$$
We notice that for a smooth open set $M'$ in $M$ the identity map is a local model of degree $1$ for $M'$ so the intersection theory on $M'$ is the usual one. 

This theory generalizes the usual intersection theory for integral analytic cycles in complex manifolds.
\item
Let $q : \tilde{M} \rightarrow M$ be a local model of degree $k$ and suppose that we have two analytic families of integral cycles $(X_{s})_{s \in S}$,  $(Y_{s})_{s \in S}$ in $M$ such that $X_{s}$ and $Y_{s}$ intersect properly in $M$ for all $s$ in $S$. Then $(q^{*}(X_{s}))_{s \in S}$ and $(q^{*}(Y_{s}))_{s \in S}$ are analytic families of cycles in $\tilde{M}$ and consequently so is $(q^{*}(X_{s}) \cap q^{*}(Y_{s}))_{s \in S}$ (see [B-M.2]). It follows that the rational cycles $X_{s}\cap_{M} Y_{s} = \frac{1}{k}.q_{*}(q^{*}(X_{s}) \cap q^{*}(Y_{s}))$ form  an analytic family of cycles in $M$. 
\hfill{$\square$}
\end{enumerate}

It is worthwhile to point out that the proof of the theorem \ref{intersect. 1} contains  the following result.

\begin{prop}\label{moving}
Let $X$ and $Y$ two integral cycles that intersect properly  in a nearly smooth complex space $M$. For any point $x_{0}\in M$ there exists a local model $q : \tilde{M'} \to M'$ on $M$ with $x_{0}\in M'$  and  two analytic families of integral cycles $(X_{s})_{s \in S}$ and $(Y_{s})_{s \in S}$ in $M'$ parametrized  by an open neighbourhood $S$ of the origin in a numerical space such that the following conditions hold:
\begin{enumerate}[i)]
\item 
We have $X_{0} = k.X \cap M'$ and $Y_{0} = k.Y \cap M'$ where $k := \deg q$.
\item 
On an open dense subset $S'$ in $S$, the generic point of $|X_{s}|\cap |Y_{s}|$ is not a branch point for $q$.
\item  
Generically on $|X_s|\cap |Y_{s}|$ the analytic subsets $|X_{s}|$ and $|Y_{s}|$ are smooth and transversal in $M'$.
\hfill{$\blacksquare$}
\end{enumerate}
\end{prop}

\subsection{The fundamental class of an intersection}

In general  trace morphisms do not commute with exterior product of differential forms, but all the same we have the following important result.

\begin{prop}\label{fund. class inter.}
Let $X$ and $Y$ be two analytic cycles that intersect properly in a nearly smooth complex space $M$. Then the fundamental classes of $X$, $Y$ and $X \cap_{M} Y$ in $M$ satisfy the following identity 
\begin{equation*}
c_{X\cap Y}^{M} \ = \ c_{X}^M\cup c_{Y}^M.
\end{equation*}
\end{prop}

\parag{Proof} 
It is enough to prove the result in the case of a local model $q : \tilde{M} \rightarrow M$,  whose degree will be denoted $k$.

\smallskip
First we will show that  
\begin{equation}\label{eq.4}
q^{*}(X) \cap_{\tilde{M}} q^{*}(Y) = q^{*}(X\cap_{M} Y)
\end{equation} 
in $\tilde{M}$, and without loss of generality we may assume $X$ and $Y$ irreducible. It is clear that equation (\ref{eq.4}) is valid if the generic point $y$ in $X\cap Y$ is not contained in the branch locus of $q$ and $X$ and $Y$ are both smooth and transversal at $y$. The general case can then be reduced to this case by the same kind of \lq moving technique\rq\ as was used in the proof of theorem \ref{intersect. 1}.

\smallskip
Now $\tilde{M}$ being smooth the identity $c_{q^*X}^{\tilde M}\cup c_{q^*Y}^{\tilde M}  = c_{q^*X\cap q^*Y}^{\tilde M}$ is valid. Hence  from  remark ii) following proposition \ref{N-Smooth 0}, theorem \ref{pull-back fund. class} and equation (\ref{eq.4}) above we get
$$
\hat{q}^*(c_{X}^M\cup c_{Y}^M) = \hat{q}^*(c_{X}^M)\cup \hat{q}^*(c_{Y}^M) = c_{q^*X}^{\tilde M}\cup c_{q^*Y}^{\tilde M} 
= c_{q^*X\cap q^*Y}^{\tilde M} = c_{q^*(X\cap Y)}^{\tilde M} = \hat{q}^*(c_{X\cap Y}^{M})
$$
and consequently we obtain $c_{X\cap Y}^{M} \ = \ c_{X}^M\cup c_{Y}^M$ since $\hat{q}^*$ is injective.
\hfill{$\blacksquare$}

\parag{Example} Put $M := \{(x, y, z) \in \C^{3} \ / \  x.y = z^{2} \}$. 
It is easy to see that $M$ is the image of the  holomorphic map 
$$
\C^{2} \longrightarrow\C^{3}\qquad (u, v) \mapsto (u^{2}, v^{2}, u.v)
$$ 
and we let $q : \C^2\rightarrow M$ denote the induced map. The map $q$ is proper, finite and surjective, so $M$ is nearly smooth and in fact $M\simeq \C^2/\{{\pm 1}\}$.

Consider the reduced cycles $X := \{ x = z = 0 \} $ and $Y := \{y = z = 0 \}$ in $M$. They are both of pure codimension one and they  intersect properly in $M$. We see that $q^{*}(X)$ is the reduced line $D := \{ u = 0\}$ in $\C^{2}$ and that $q^{*}(Y)$ is the reduced line $D' := \{ v = 0 \}$  in $\C^{2}$. These two lines  are transversal at the origin and consequently $D \cap_{\C^2} D' = 1.\{0\}$. As $q_{*}(1.\{0\}) = 1. \{0\}$ (because $q$ is bijective from $\{0\}$ to $\{0\}$) we obtain $X \cap_{M} Y = \frac{1}{2}.\{0\}$.

\smallskip
Let us now compute the fundamental classes (in Cech cohomology).  From the fact that $\frac{1}{2}.Trace_{q}[q^*\alpha] = \alpha$, for every holomorphic form on an open set in $M$, we deduce the following identities (for classes of Cech cohomology):
$$
Trace_{q}\left[\frac{du}{u}\right] = \left[\frac{dx}{x}\right], \ Trace_{q}\left[\frac{dv}{v}\right] = \left[\frac{dy}{y}\right] \  \text{and} \ \ 
Trace_{q}\left[\frac{du}{u}\wedge \frac{dv}{v}\right] = \frac{1}{2}\left[\frac{dx}{x}\wedge \frac{dy}{y}\right].
$$
Hence we get
$$
c_{X}^{M} \ = \ \frac{1}{2}.Trace_{q}[c_{q^*X}^{\C^2}] \ = \ \frac{1}{2}.Trace_{q}\left[\frac{1}{2i\pi}\frac{du}{u}\right] \ = \ \frac{1}{2}\left[\frac{1}{2i\pi}\frac{dx}{x}\right]
$$
and in the same way $c_{Y}^{M} = \frac{1}{2}\left[\frac{1}{2i\pi}\frac{dy}{y}\right]$. It then follows that 
$$
c_{X}^{M} \cup c_{Y}^{M} \ = \ \frac{1}{4}\left[\frac{1}{2i\pi}\frac{dx}{x}\right]\cup\left[\frac{1}{2i\pi}\frac{dy}{y}\right] \ = \ \frac{1}{4}\left[\frac{1}{(2i\pi)^2}\frac{dx}{x}\wedge\frac{dy}{y}\right].
$$
On the other hand we obtain
\begin{eqnarray*}
c_{X \cap Y}^{M} = c_{\frac 12.\{0\}}^{M} &=& \frac 12.c_{\{0\}}^{M} = \frac{1}{4}.Trace_{q}\left[c_{q^*\{0\}}^{\C^2}\right] = \frac{1}{4}.Trace_{q}\left[c_{2\{0\}}^{\C^2}\right]\\ 
&=& \frac 12.Trace_{q}\left[\frac{1}{(2i\pi)^2}\frac{du}{u}\wedge \frac{dv}{v}\right] = 
\frac{1}{4}\left[\frac{1}{(2i\pi)^2}\frac{dx}{x}\wedge\frac{dy}{y}\right]
\end{eqnarray*}
so indeed we have  $c_{X}^{M} \cup c_{Y}^{M} =  c_{X \cap Y}^{M}$.

\parag{Remark}
Since intersection of analytic cycles in a complex manifold is both commutative and associative, we easily deduce from formula (\ref{eq.4}) that the same is true for intersection of analytic cycles in a nearly smooth complex space. 

\subsection{Generalization of $X._{f}Y$}

\begin{thm}\label{$X._{f}Y$}
Let $f : M \rightarrow N$ be a holomorphic map between complex spaces where $M$ is of pure dimension and $N$ is nearly smooth.  Let $X$ be a cycle  in $M$ and $Y$ be a cycle in $N$ such that $f^{-1}(|Y|) \cap |X| = \emptyset$ or
$$
\codim_Mf^{-1}(|Y|) \cap |X| = \codim_{M}X + \codim_{N}Y.
$$
Then there exists a unique cycle in $M$, denoted $X._fY$, characterized by the following properties:
\begin{enumerate}[i)]
\item
$|X._fY| = f^{-1}(|Y|) \cap |X|$
\item
For every open subset $M'$ of $M$ one has
$$
(X._fY)\cap M'  = (X\cap M')._{f_{|M'}}Y
$$
\item
For every closed embedding $j : M'\rightarrow M_1$ from an open subset of $M$ into a complex manifold one has 
$$
(j,f_{|M'})_*((X._fY)\cap M') = (j,f_{|M'})_*(X\cap M')\cap_{M_1\times N}(M_1\times Y).
$$
\end{enumerate}
\end{thm}

\parag{Proof} In the situation of condition iii) the map $(j,f_{|M'}) : M' \rightarrow M_1\times N$ is an embedding and if $f^{-1}(|Y|) \cap |X| \not= \emptyset$ we have
\begin{align*}
\codim_{M_{1}\times N} (j,f_{|M'})(|X|) \cap &(M_{1}\times |Y|) = \\
&\codim_{M_{1}\times N} (j,f_{|M'})(|X|) + \codim_{M_{1}\times N} (M_{1}\times |Y|).
\end{align*}
It follows that the right hand side of the equation in condition iii) is well defined and determines a cycle in $M'$. Hence to prove the theorem it is sufficient to show that this cycle does not depend on the choice of the embedding $j$ because we then get the cycle $X._fY$ by the following proceedure. We cover $M$ with open subsets each of which is endowed with a closed embedding into some complex manifold. Then the family of cycles in these subsets will automatically glue together and give us the cycle $X._fY$.

Now the proof of the independence of $j$ reduces to showing that, for a closed embedding  
$J : M_1\rightarrow M_2$ of one manifold into another and a cycle $Z$ in $M_1\times N$, we have
\begin{equation}\label{eq.5}
(J\times\id_N)_*(Z\cap_{M_1\times N}(M_1\times Y)) = (J\times\id_N)_*(Z)\cap_{M_2\times N}(M_2\times Y).
\end{equation} 
In order to show this we first remark that (\ref{eq.5}) is valid in the case where $N$ is smooth (see [B-M.2]). Then we notice also that the assertion is local on $N$ so we may assume that we have a local model $q :  \tilde{N} \rightarrow N$ of a certain degree $k$.
Now by considering the commutative diagram
$$ 
\xymatrix{
M_{1} \times \tilde{N} \ar[rr]^{J\times \id_{\tilde{N}}} \ar[d]_{\id_{M_1}\times q} && M_{2} \times\tilde{N} \ar[d]^{\id_{M_2}\times q}\\
M_{1}\times N \ar[rr]^{J\times\id_N} && M_{2} \times N} 
$$
observing that $(J\times \id_{\tilde{N}})_{*}(\id_{M_1}\times q)^*Z = (\id_{M_2}\times q)^*(J\times\id_N)_*Z$ and keeping (\ref{eq.4}) in mind we get
\begin{align*}
(J\times\id_N)_*&(Z\cap_{M_1\times N}(M_1\times Y)) \\ 
&= \frac 1k(J\times\id_N)_*(\id_{M_1}\times q)_*(\id_{M_1}\times q)^*(Z\cap_{M_1\times N}(M_1\times Y))\\
&= \frac 1k(\id_{M_2}\times q)_*(J\times\id_{\tilde N})_*[(\id_{M_1}\times q)^*(Z)\cap_{M_1\times N}(M_1\times q^*Y)]\\
&= \frac 1k(\id_{M_2}\times q)_*[(J\times\id_{\tilde N})_*(\id_{M_1}\times q)^*(Z)\cap_{M_1\times N}(M_2\times q^*Y)]\\
&= \frac 1k(\id_{M_2}\times q)_*[(\id_{M_2}\times q)^*(J\times\id_N)_*(Z)\cap_{M_1\times N}(\id_{M_2}\times q)^*(M_2\times Y)]\\
&= (J\times\id_N)_*(Z)\cap_{M_2\times N}(M_2\times Y).
\end{align*}
This completes the proof.
\hfill{$\blacksquare$}\\

The cycle $X._fY$ is called the {\bf $f-$intersection product} of $X$ and $Y$.

\parag{Example} Let $f : M \rightarrow N$ be a geometrically flat holomorphic map from a complex space to a nearly smooth complex space\footnote{Note that, as $N$ is normal, this just means that $f$ is  equidimensional.}. Then every fibre of $f$ is of dimension $\dim M - \dim N$ so for any cycle $Y$ in $N$ we have 
$$
\dim f^{-1}(\vert Y\vert) = \dim M - \dim N + \dim Y
$$ 
and consequently the cycle $M._{f}Y$ is defined. We shall show in lemma \ref{racc.} that in this case we have $M._{f} Y = f^{*}(Y)$.
\hfill$\square$\\

\begin{thm}\label{parameters}
Let $f : M \rightarrow N$ be a holomorphic map from a pure dimensional complex space to a nearly smooth space. Let $(X_{s})_{s \in S}$ and $(Y_{s})_{s \in S}$ be analytic families of cycles in $M$ and $N$ parametrized by a reduced complex space $S$. Assumme that for each $s \in S$ the cycle $X_{s}._{f}Y_{s}$ is well defined in $M$. Then the family of cycles $(X_{s}._{f}Y_{s})_{s \in S}$ is analytic.
\end{thm}
\parag{Proof} We may assume that we have a closed  embedding $j : M\rightarrow M_1$ of $M$ into a complex manifold. Then from [B-M.2] we know that $((j,f)_*(X_s))_{s\in S}$ and $(M_1\times Y_s))_{s\in S}$ are analytic families of cycles in the nearly smooth space $M_1\times N$, and remark ii) following definition \ref{intersect.2} allows us to conclude.
$\hfill \blacksquare$

\begin{lemma}\label{racc.}
Let $f : M \rightarrow N$ be a geometrically flat holomorphic map from a pure dimensional complex space $M$ to a nearly smooth complex space $N$. Then for any cycle $Y$ in $N$ we have
$$
f^{*}(Y) = M._{f}Y.
$$
\end{lemma}

\parag{Proof} Recall that, by definition, $f^{*}(Y)$ is the  graph cycle of the analytic family of fibres of $f$ parametrized by $Y$ when $Y$ is reduced.\\
 Both sides being additive in $Y$ it is enough to consider the case where $Y$ is an irreducible subset in $N$. From [B-M.2] we know that the result is valid in the case where $N$ is smooth and since $f$ is equidimensional it is also valid in the case where $Y$ is not contained in the singular part, $S(N)$, of $N$. The problem is local in $N$ so, by proposition \ref{moving}, we may assume that there exists an analytic family of cycles $(Y_s)_{s\in D}$ in $N$, where $D$ is an open disc centered at the origin in $\C$, such that for all $s$ in an open  dense subset $D'$ of $D$ no irreducible component of $|Y_s|$ is contained in $S(N)$  and $Y_0 = Y$. Then $f^*(Y_s) = M._fY_s$ for all $s\in D'$ and consequently $f^*(Y) = M._fY$ by continuity, thanks to the previous theorem.
$\hfill \blacksquare$

\begin{thm}[projection formula]\label{proj. form.}
Let $f : M \to N$ be a holomorphic map between a complex space $M$ and a nearly smooth complex space $N$. Let $X$ a cycle in $M$ and $Y$ a cycle in $N$ such that the cycle $X._{f}Y$ is well defined. Assume also that the restriction of $f$ to the support of $X$ is proper. Then the cycles $f_{*}(X)$ and $Y$ are intersecting properly in $N$ and  the projection formula  $ f_{*}(X._{f}Y) = f_{*}(X)\cap_{N}Y $ holds.
\end{thm}

\parag{Proof}
As our problem is local on $N$, we may assume that we have a local model $q : \tilde{N} \rightarrow N$ of degree $k$. 
But as the problem is also local on $M$, we may assume that we have a closed embedding $j : M \hookrightarrow M_{1}$ of $M$ into a complex manifold $M_{1}$.

Consider the following diagram
$$ 
\xymatrix{\quad & M_{1} \times \tilde{N} \ar[r]^{\quad\tilde{p}} \ar[d]^{\id\times q} & \tilde{N} \ar[d]^{q}\\
M \ar[r]^{(j,f) \quad } & M_{1}\times N \ar[r]^{\quad p} & N} 
$$
where $p$ and $\tilde{p}$ are the natural maps. The identity $p\circ (j\times f) = f$ is obvious and for any cycle $Z$ in $M_{1}\times N$ such $p_{*}(Z)$ is defined we have the equality
\begin{equation}\label{eq.6}
\tilde{p}_{*}(id\times q)^{*}(Z) = q^{*}(p_{*}(Z)). 
\end{equation}
By additivity, it is enough to prove (\ref{eq.6}) for an irreducible $Z$. In this case if $d$ is the generic degree of the induced map $Z \rightarrow p(Z)$, $q^{*}(p_{*}(Z))$ is $d$ times the graph cycle of the family of fibres of $q$ parametrized by the set $\vert p_{*}(Z)\vert$, and the cycle $\tilde{p}_{*}\big((\id\times q)^{*}(Z)\big)$ is clearly equal to $p_{*}(Z)$, as $(\id\times q)^{*}(Z)$ is the graph cycle of the family of the fibres of $(\id\times q)$ parametrized by $Z$.\\

Note that $\id\times q : M_{1}\times \tilde{N} \to M_{1}\times N$ is a local model of degree $k$  for the nearly smooth space $M_{1}\times N$. Then we have, with $Z := (j\times f)_{*}(X)$:
\begin{align*}
& f_{*}(X._{f}Y) = p_{*}\big((j,f)_{*}(X._{f}Y)\big) = p_{*}\Big((j,f)_{*}(X)\cap_{M_{1}\times N}p^{*}(Y)\Big) \\
& \quad = p_{*}\Big(\frac{1}{k}(\id\times q)_{*}\big((\id\times q)^{*} (Z)\cap_{M_{1}\times \tilde{N}} (\id\times q)^{*}(p^{*}(Y) \big)\Big)\\
& \quad  = \frac{1}{k} q_{*}\Big(\tilde{p}_{*}\big((\id\times q)^{*}(Z) \cap_{M_{1}\times \tilde{N}} (\id\times q)^{*}(p^{*}(Y)\big)\Big) \quad {\rm using} \quad  p\circ (\id\times q) = q \circ \tilde{p} \\
& \quad =   \frac{1}{k} q_{*}\Big(\tilde{p}_{*}\big((\id\times q)^{*}(Z)  \cap_{\tilde{N}} \tilde{p}^{*}(q^{*}(Y)\big)\Big) ; \ {\rm  \ the \ projection \ formula \ for} \ \tilde{p} \ {\rm gives:}\\
& \quad =  \frac{1}{k} q_{*}\Big( \tilde{p}_{*}\big((\id\times q)^{*}(Z)\big) \cap_{\tilde{N}} q^{*}(Y)\Big) \quad {\rm and \ by} \quad (\ref{eq.6})\\
& \quad =  \frac{1}{k} q_{*}\Big(q^{*}(p_{*}(Z)) \cap_{\tilde{N}}q^{*}(Y)\Big) = p_{*}(Z) \cap_{N} Y = f_{*}(X)\cap_{N}Y.
\end{align*}
This concludes the proof.
 $\hfill \blacksquare$\\
 
 The next theorem shows an important compatibility property between the pull-back of cycles, $f-$intersection and the intersection in a nearly smooth complex space, generalizing the standard analogous  compatibility for complex manifolds.

 \begin{thm}\label{compat.1}
Let $f : M \to N$ be a holomorphic map between two nearly smooth complex spaces. Let $X$ be a cycle in $M$ and $Y$ a cycle in $N$ such that the cycles $M._fY$ and  $X._{f}Y$ are defined. Then the cycles $M._fY$ and  $X$ intersect properly in $M$ and
\begin{equation}\label{eq.7}
 X._{f}Y \ = \ X \cap_{M}(M._fY). 
 \end{equation}
\end{thm}

\parag{Proof} The cycle $X._{f}Y$ is defined when
$$
\codim_{M} X + \codim_{N} Y = \codim_{M} (|X|\cap f^{-1}(|Y|)) 
$$
in a neighbourhood of $|X|\cap f^{-1}(|Y|)$ and the cycle $M._fY$ is defined when $ \codim_{M}  f^{-1}(|Y|)$ is equal to $\codim_{N} Y $. Hence, up to replace $M$ by an open neighbourhood of $|X|\cap f^{-1}(|Y|)$ we get either   $|X|\cap f^{-1}(|Y|) = \emptyset$  or 
$$
\codim_{M} X + \codim_{M}  f^{-1}(\vert Y\vert) = \codim_{M} (|X|\cap f^{-1}(|Y|)) 
$$
so the cycles $M._fY$ and $X$ intersect properly in $M$ and the first assertion of the theorem is proved.\footnote{As we assume that $M$ is nearly smooth, we have 
$$
\codim_{M} (|X|\cap f^{-1}(|Y|)) \leq \codim_{M} X + \codim_{M}  f^{-1}(|Y|)
$$
and as  $\codim_{M}  f^{-1}(|Y|) \leq \codim_{N} Y$ the hypothesis that $M._fY$ is defined near $|X|\cap f^{-1}(|Y|)$ is a consequence of the hypothesis that $X._{f}Y$ is defined and of the fact that $M$ is nearly smooth.}

To prove the identity (\ref{eq.7}) we may assume that we have an embedding $j : M\rightarrow M_1$ where $M_1$ is a complex manifold and in the sequal we will refer to  the following commutative diagram
$$
\xymatrix{
M\ar[rr]^{(\id_M,f)\quad}\ar[drr]_{(j,f)}&&M \times N \ar[rr]^{\pi} \ar[d]^{j\times \id_{N}} && M\ar[d]^{j}\\
&&M_{1}\times N \ar[rr]^{\pi_1} && M_{1}} 
$$
where $\pi$ and $\pi_1$ are the natural projections. Notice that $(\id_M,f)_*M$ is the graph of the map $f$ and will be denoted $G_f$.

\begin{lemma}\label{compat.0}
In the situation above we have the following identities:
\begin{enumerate}[i)]
\item
\ $(\id_M,f)_*X \ = \ (X\times N)\cap_{M\times N} G_f$.
\item 
\ $\pi_*\left[(\id_M,f)_*X\cap_{M\times N}(M\times Y)\right] \ = \ X._fY$.
\end{enumerate}
\end{lemma}
\parag{Proof} To prove i) we use the projection formula for $\pi$ to get
$$
\pi_*((X\times N)\cap_{M\times N} G_f) \ = \ \pi_*(\pi^*X\cap_{M\times N} G_f) \ = \ X\cap_M\pi_*G_f = X\cap_MM \ = \ X
$$
and consequently $(\id_M,f)_*X \ = \ (X\times N)\cap_{M\times N} G_f$ since $\pi\circ(\id_M,f) = \id_M$.

\smallskip
To prove ii) we use the commutativity of the diagram above and the projection formula for $j\times\id_N$ to get
\begin{align*}
j_*\pi_*\left[(\id_M,f)_*X\cap_{M\times N}(M\times Y)\right] = (\pi_1)_*(j\times\id_N)_*\left[(\id_M,f)_*X\cap_{M\times N}(M\times Y)\right]&\\
= (\pi_1)_*(j\times\id_N)_*\left[(\id_M,f)_*X\cap_{M\times N}((M\times Y)._{j\times\id_N}(M_1\times Y))\right]&\\
= (\pi_1)_*\left[(j\times\id_N)_*(\id_M,f)_*X\cap_{M\times N}(M_1\times Y)\right]&\\
= (\pi_1)_*\left[(j,f)_*X\cap_{M\times N}(M_1\times Y)\right] = j_*(X._fY)&
\end{align*}
and hence $\pi_*\left[(\id_M,f)_*X\cap_{M\times N}(M\times Y)\right] \ = \ X._fY$.
\hfill{$\blacksquare$}\\

We now deduce (\ref{eq.7}) from lemma \ref{compat.0} by using the associativity of intersection and the projection formula for $\pi$:
\begin{eqnarray*}
X._{f}Y &=& \pi_*[(\id_M,f)_*X\cap_{M\times N}(M\times Y)]\\
&=& \pi_*[((X\times N)\cap_{M\times N}G_f)\cap_{M\times N}(M\times Y)]\\
&=& \pi_*[(X\times N)\cap_{M\times N}(G_f\cap_{M\times N}(M\times Y))]\\ 
&=& \pi_*[\pi^*(X)\cap_{M\times N}(G_f\cap_{M\times N}(M\times Y))]\\ 
&=& X\cap_{M}\pi_*[G_f\cap_{M\times N}(M\times Y)]\\ 
&=& X \cap_{M}(M._fY). 
\end{eqnarray*}
\hfill $\blacksquare$\\
 
The following corollary is an interesting application of the previous theorem and the projection formula.

\begin{cor}\label{le point}
Let $M$ be a nearly smooth complex space, $P$ an analytic subset of $M$ and suppose that $P$ is nearly smooth. Let $X$ be a cycle in $P$ and $Y$ a cycle in $M$ such that $Y$ intersects both $X$ and $P$  properly in $M$. Then the cycle $X$ and $P \cap_{M} Y$ intersect properly in $P$ and we have
\begin{equation}\label{eq.8} 
X\cap_{M} Y = X \cap_{P}(P \cap_{M}Y).
\end{equation}
 \end{cor}
 
\parag{Proof} It is clear that $X$ and $P \cap_{M} Y$ intersect properly in $P$.

\smallskip
If we let $j: P \hookrightarrow M$ denote the inclusion map, then the precise meaning of (\ref{eq.8}) is
$$
j_*X\cap_{M} Y = j_*[X \cap_{P}(P._jY)].
$$
But from theorem \ref{compat.1} and the projection formula for $j$ we get
$$
j_*[X \cap_{P}(P._jY)] = j_*[X._jY)] = j_*X\cap_{M} Y
$$
and that completes the proof.
\hfill{$\blacksquare$}\\

We only sketch the proof of the following theorem leaving standard details to the reader. 
\begin{thm}\label{assoc.+}
Assume that we have holomorphic maps $M \overset{f}{\to} N \overset{g}{\to} P $ where $M$ is  a reduced complex space and where  $N$ and $P$ are nearly smooth complex spaces. For cycles $X, Y, Z$ respectively in $M, N, P$ such that $X._{f}Y, Y._{g}Z$ and $X._{f}(Y._{g}Z)$ are defined, the cycle $(X._{f}Y)._{g\circ f}Z$ is defined and we have
$$ (X._{f}Y)._{g\circ f}Z = X._{f}(Y._{g}Z).$$
\end{thm}

\parag{Sketch of proof} Let $j : M \to M_{1}$ be a local closed embedding of   $M$ in a complex manifold  $M_{1}$. The result is consequence of the associativity of the intersection in the nearly smooth complex space  $M_{1}\times N\times P$ for the three cycles  $ (j, f, g\circ f)_{*}X, \  (\id_{M_{1}\times N}, \ g)_{*}(M_{1}\times Y) \ {\rm and}  \  M_{1}\times N\times Z$ which meet properly in  $M_{1}\times N\times P$,  noticing that the intersections 
\begin{itemize}
\item  $  (\id_{M_{1}\times N}, g)_{*}(M_{1}\times Y) \cap (M_{1}\times N \times Z) $
\item $ \big( (j, f, g\circ f)_{*}X\big) \cap  \big(\id_{M_{1}\times N}, g)_{*}(M_{1}\times Y) \cap (M_{1}\times N \times Z)\big) $
\item $   (j, f, g\circ f)_{*}X \cap  (\id_{M_{1}\times N}, g)_{*}(M_{1}\times Y) $
\item $  \big((j, f, g\circ f)_{*}X \cap  (\id_{M_{1}\times N}, g)_{*}(M_{1}\times Y)\big) \cap (M_{1}\times N \times Z)$
\end{itemize}
correspond respectively to the cycles 

$$ M_{1}\times (Y._{g}Z), \quad   X._{f}(Y._{g}Z), \quad    (j, f, g)_{*}(X._{f}Y)\quad  {\rm and}  \quad  (X._{f}Y)._{g \circ f} Z. $$

\subsection{Some formulae for direct image of fundamental classes}\

As in the case of a holomorphic map between complex manifolds, the computation of the fundamental classe of a direct image cycle will be consequence of two special cases : the case of a (local) embedding in a nearly smooth complex space and the case of the natural projection $\pi : M_{1}\times N \to N$ where $M_{1}$ is a complex manifold and $N$ is a nearly smooth complex space.
 This is enough to follow the computation of  fundamental classes in the projection formula (see theorem \ref{proj. form.}).\\
 So we begin with the behavior of the fundamental class in a closed embedding of nearly smooth complex spaces.

\begin{prop}\label{plong. inv.}
Let $j : M \to N$ be a closed holomorphic embedding between two nearly smooth complex spaces $M$. Assume that $M$ and $N$ are have pure dimensions and that $dim N - dim M := r$. Then the fundamental class $c_{j_{*}(M)}^{N}$ of $j_{*}(M)$ in $N$ induces a sheaf morphism
$$ \gamma^{\bullet} : j_{*}(\omega_{M}^{\bullet}) \to \underline{H}_{j(M)}^{r}(\omega_{N}^{r+\bullet}) $$
and for each co-dimension $p$ cycle $X$ in $M$ the fundamental classes of $X$ in $M$ and the fundamental class of $j_{*}(X)$ in $N$  are related by the equality
$$ H_{j(\vert X\vert)}^{p}(\gamma)[c_{X}^{M}] = c_{j_{*}(X)}^{N},$$
where we use the natural isomorphism $\underline{H}_{j(\vert X\vert)}^{p}(\underline{H}_{j(M)}^{r}(\omega_{N}^{r+p}) ) \simeq \underline{H}_{j(\vert X\vert)}^{p+r}(\omega_{N}^{p+r})$ which is consequence of the vanishing of the sheaves $\underline{H}_{j(M)}^{s}(\omega_{N}^{\bullet})$ for $s \not= r$.
\end{prop}

\parag{Proof} It is easy to see that  $\underline{H}_{j(\vert X\vert)}^{p}(\gamma)[c_{X}^{M}] $ satisfies the conditions of the lemma \ref{caract.}, thanks to the transitivity of the topological fundamental classes;  then this lemma  gives the equality.$\hfill \blacksquare$\\

\begin{prop}\label{product case}
Let $M_{1}$ be a connected complex manifold of dimension $p_{1}$ and $N$ be a nearly smooth complex space. Let $Z$ be a codimension $p_{1}+p$ cycle in $M_{1}\times N$ such that the restriction to $\vert Z\vert$ of the projection $\pi : M_{1}\times N \to N$ is proper. Then there exists a natural morphism of direct image
$$ \theta : \pi_{*}\big(\underline{H}_{\vert Z\vert}^{p_{1}+p}(\omega_{M_{1}\times N}^{p_{1}+p}\big) \to \underline{H}_{\vert Y\vert}^{p}(\omega_{N}^{p}) $$
sending the fundamental class of $Z$ in $M_{1}\times N$ to the fundamental class  in $N$  of the cycle  $Y := \pi_{*}(Z)$.
\end{prop}

\parag{Proof} The assertion is local on $N$ so we may assume that we have a local model $q : \tilde{N} \to N$. Then $\tilde{q} := (id_{M_{1}}\times q) : M_{1}\times \tilde{N} \to M_{1}\times N$ is a local model for the nearly smooth complex space $M_{1}\times N$. The corresponding result for complex manifolds gives, for the cycle $\tilde{Z} := \tilde{q}^{*}(Z)$,   a morphism
$$\tilde{\theta} : \tilde{\pi}_{*}\big(\underline{H}_{\vert \tilde{Z}\vert}^{p_{1}+p}(\Omega_{M_{1}\times \tilde{N}}^{p_{1}+p})\big) \to \underline{H}_{\vert \tilde{Y}\vert}^{p}(\Omega_{\tilde{N}}^{p}) $$
where we note $\tilde{\pi}$ the projection $M_{1}\times \tilde{N} \to \tilde{N}$ and  $\tilde{Y} := \tilde{\pi}_{*}(\tilde{Z})$. This morphism induces the direct image of $\bar\partial-$closed currents with supports in $\vert \tilde{Z}\vert$ (resp. in $\vert \tilde{Y}\vert$)  in the corresponding moderate cohomologies and  sends the integration current  of $\tilde{Z}$ in $M_{1}\times \tilde{N}$ to the integration current  in $\tilde{N}$  of the cycle  $\tilde{Y}$. It is then an easy  to conclude using the direct image by $q$ and  proposition \ref{N-Smooth 0}.$\hfill \blacksquare$\\

\parag{References}
\begin{itemize}
\item{[B.75]} Barlet, D. \textit{Espace analytique r\'eduit des cycles analytiques compacts d'un espace analytique complexe de dimension finie},  in Sem. F. Norguet II, Lecture Notes, vol. 482 Springer Verlag (1975), pp. 1-158.
\item{[B.78]} Barlet, D. {\it Le faisceau $\omega_{X}^{\bullet}$ sur un espace complexe r\'eduit de dimension pure},  in Sem. F. Norguet III, Lecture Notes, vol. 670, Springer Verlag (1978), p. 187-204.
\item{[B.80]} Barlet, D. \textit{Familles analytiques de cycles et classes fondamentales relatives},  in Sem. F. Norguet IV, Lecture Notes, vol. 807 Springer Verlag (1980), pp. 1-24.
\item{ [B-M 1]} Barlet, D. et Magnusson, J. {\it Cycles analytiques complexes I: th\'eor\`{e}mes de pr\'eparation des cycles}, Cours sp\'ecialis\'es 22, SMF 2014.
\item{ [B-M 2] }Barlet, D. et Magnusson, J. {\it Cycles analytiques complexes II},  to appear.
\item{[D. 69]} Draper, R. {\it Intersection theory in complex analytic geometry},  Math. Ann. 180 (1969) p. 175-204.
\item{[Fult.]} Fulton, W. {\it Intersection theory}. Second edition. Springer-Verlag, Berlin, 1998.  
\item{[F-MP]}  Fulton, W. and  MacPherson, R. {\it Defining algebraic intersections}. Algebraic geometry (Proc. Sympos., Univ. Troms¿, Troms¿, 1977), pp.1-30, Lecture Notes in Math., 687, Springer, Berlin, 1978.
\item{[K.09 ]} Koll\'{a}r, J. {\it Lectures on Resolution of Singularities}, Priceton University press (2009).
\item{[T.95]} Tworzewski, P. {\it Intersection theory in complex analytic geometry},  Ann. Polon. Math. 62 (1995), no. 2, p. 177-191.

\end{itemize}

\end{document}